\documentclass[12pt]{article}
\usepackage{mathrsfs,amsmath,amssymb,amsfonts}
\numberwithin{equation}{section}

\usepackage[latin1]{inputenc}
\usepackage{multirow}
\usepackage[justification=centering]{caption}
\usepackage{amsfonts}
\usepackage{float}
\usepackage{subfigure}
\usepackage{stmaryrd}
\usepackage{amsmath}
\allowdisplaybreaks[4]
\usepackage{booktabs}
\usepackage{adjustbox}
\usepackage{multicol}
\usepackage{bm}
\usepackage{upgreek}
 \usepackage{graphicx} %use graph format
  \usepackage{epstopdf}
\usepackage[usenames]{color}
\usepackage{caption}

\usepackage{subfigure}
\usepackage{enumerate}
\usepackage{amsthm}
\newcommand{\lbl}{\label}

\newcommand{\ignore}[1]{}{}
\newcommand{\be}{\begin{equation}}               %\be=\begin{equation}
\newcommand{\ee}{\end{equation}}                 %\ee=\end{equation}
\newcommand{\bi}{\begin{itemize}}
\newcommand{\ei}{\end{itemize}}

\newcommand{\beaa}{\begin{eqnarray*}}
\newcommand{\eeaa}{\end{eqnarray*}}
\newcommand{\bea}{\begin{eqnarray}}
\newcommand{\eea}{\end{eqnarray}}

\newcommand{\beqn}{\begin{eqnarray}}             %\beqn=\begin{eqnarray}
\newcommand{\eeqn}{\end{eqnarray}}               %\eeqn=\end{eqnarray}
\newcommand{\beq}{\begin{eqnarray*}}             %\beq=\begin{eqnarray*}
\newcommand{\eeq}{\end{eqnarray*}}               %\eeq=\end{eqnarray*}
\newcommand{\ssb}{\scriptstyle \footnotesize % \scriptsize
                 \begin{array}{c}}
\newcommand{\esb}{\end{array}}
\newtheorem{theorem}{{ Theorem}}%[section]
 \newtheorem{lemma}{{ Lemma}}%[section]
 \newtheorem{remark}{{Remark}}%[section]
\makeatletter

\newcommand{\Rmnum}[1]{\expandafter\@slowromancap\romannumeral #1@}
\makeatother
\begin{document}
\title{Variable selection for partially linear single-index varying-coefficient model
\footnote{
Liugen Xue's work is supported by the National Natural Science Foundation of China under Grant [No. 12471252];
Junshan Xie's work is supported by  the Natural Science Foundation of Henan Province of China under Grant [No. 242300421379].
}
 }
\author{ Lijuan Han, Liugen Xue, Junshan Xie\thanks{\textsl{E-mail address}:
\texttt{junshan@henu.edu.cn.}}}
\date{\small \it School of Mathematics and Statistics, Henan University,
Kaifeng, 475000, P.R. China
} \maketitle
\mbox{}\hrule\mbox{}\\[0.5cm]
\textbf{Abstract}\\[-0.2cm]^^L
%In this paper, we
This paper focuses on variable selection for a partially linear single-index varying-coefficient model. A regularized variable selection procedure by combining basis function approximations with SCAD penalty is proposed. %The proposed procedure
 It can simultaneously select significant variables in the parametric %components
  and nonparametric components and estimate the nonzero regression coefficients and coefficient functions. %We establish the
 The consistency of the variable selection procedure and the oracle property of the penalized least-squares estimators for high-dimensional data  are established. %by choosing the appropriate tuning parameters.
   Some simulations and the real data analysis are constructed to illustrate
  the finite sample performances of the proposed method.

%---------------------------- -----------------------------------------------------------------------------
[0.5cm]
\mbox{}\hrule\mbox{}
%------------------------------------------------------------------------------------------
\newpage
%------------------------------------------------------------------------------------------

\section{Introduction}
\hspace{2em}Consider the partially linear single-index varying-coefficient model
\begin{equation}\label{1.1}
Y_{i}=\theta^{T}U_{i}+\boldsymbol{g}^{T}(\beta^{T}X_{i})Z_{i}+\varepsilon_{i},\quad i=1,\:\ldots,n,
\end{equation}
where $U_i=(U_{i1},\ldots,U_{id})^T,$ $X_i=(X_{i1},\ldots,X_{ip})^T,$ and $Z_i=(Z_{i1},\ldots,Z_{iq})^T$ are covariates. $Y_i$'s are the response variables, $\theta$ and $\beta$ are $d\times1$ and $p\times1$ vectors of unknown parameters, $\boldsymbol{g}(\cdot)$ is a $q\times1$ vector of unknown coefficient functions and the error variables $\{\varepsilon_i,~1\leq i\leq n\}$ are independent of $\{(U_i, X_i, Z_i),~ 1\leq i\leq n\}$ with the conditions that $E(\varepsilon_i)=0$, $Var(\varepsilon_i)=\sigma^2<\infty$. Generally, the first component of Z$_i$ may be taken as $1$ so that the model has an intercept function term. For the identifiability of model (1.1), it is often assumed that $\|\beta\|=1$ and the first nonzero component of $\beta$ is positive, where $\|\cdot\|$ denotes the Euclidean metric.

As an important semiparametric model, the partially linear single-index varying-coefficient model includes many other major statistical models.
If $\theta=0$, model (\ref{1.1}) reduces to the single-index varying-coefficient model studied by many scholars such as Xia and Li (1999) \cite{Xiayc}, Wu et al. (2011) \cite{Wu}, Xue and Wang (2012) \cite{Xue1}, Xue and Pang (2013) \cite{Xue2} and Wang et al. (2022) \cite{Wangt}.
If $Z_i=1$, model (\ref{1.1}) reduces to the partially linear single-index model, see Carroll et al. (1997) \cite{Carroll}, Yu and Ruppert (2002) \cite{Yu}, Xia and H$\rm{\ddot{a}}$rdle (2006) \cite{Xiac}, Wang et al. (2010) \cite{Wangj}, Liang et al. (2010) \cite{Liang}, Li et al. (2015) \cite{Lig}, Lin et al. (2022) \cite{Linh}.
If $\theta=0$ and $Z_i=1$, model (\ref{1.1}) reduces to the single-index model considered by H$\rm{\ddot{a}}$rdle et al. (1993) \cite{Hardle}, Xia et al. (2004) \cite{Xia}, Xue and Zhu (2006) \cite{Xue4}, Zhu et al. (2011) \cite{Zhu}, Peng and Huang (2011) \cite{Peng}, Li et al. (2014) \cite{Li}, Liu et al. (2019) \cite{Liu}, He et al. (2025) \cite{He}, among others.
If $\beta=1$, model (\ref{1.1}) reduces to the partially linear varying-coefficient model investigated by many statistical researchers, for example, Ahmad et al. (2005) \cite{Ahmad}, Fan and Huang (2005) \cite{Fanj}, Wang et al. (2009) \cite{Wanghj}, Zhao (2010) \cite{ZHX}, Feng and Xue (2014) \cite{Fengsy} and Zhao et al. (2023) \cite{Zhaom}.
If $\beta=1$ and $Z_i=1$, model (\ref{1.1}) reduces to the partially linear model, see Heckman (1986) \cite{Heckman}, Speckman (1988) \cite{Speckman}, Fan and Li (2004) \cite{Fanjql}, Liang and Li (2009) \cite{Liangh}, Liu et al. (2018) \cite{Liujy}, Lu et al. (2024) \cite{Luy}, etc.
If $\beta=1$ and $\theta=0$, model (\ref{1.1}) reduces to the varying-coefficient model, some works include: Hastie and Tibshirani (1993) \cite{Hastie}, Hoover et al. (1998) \cite{Hoover}, Fan and Zhang (1999) \cite{FanZ}, Huang et al. (2004) \cite{HuangW}, Xue and Zhu (2007) \cite{XueZ}, Xue and Qu (2012) \cite{XueQ},  Xiong et al. (2023) \cite{XiongT}, etc.

Model (\ref{1.1}) has three principal advantages as follows: First, it can overcome the well-known phenomenon of ``curse-of-dimensionality" which is often encountered in multivariate nonparametric settings, since $\boldsymbol{g}(\cdot)$ is a function vector of univariate variable. Second, it combines the single-index model with the varying-coefficient model so that it has a wider range of applicability. Third, it allows both discrete and continuous covariates due to the model structure, the covariates of the single-index varying-coefficient part of the model are required to be continuous, while the covariates of the linear part can be either continuous or discrete. In fact, model (\ref{1.1}) has been studied by some literature. Feng and Xue (2015) \cite{Feng} considered the problem of model detection and estimation for the single-index varying-coefficient model based on the penalized spline estimation. They used the minimum average variance estimation (MAVE) to establish the asymptotic properties of the estimators. Zhao et al. (2019) \cite{ZhaoYX} adopted a stepwise estimation procedure to estimate the index parameters, the coefficient parameters, and the coefficient functions.  Xue (2023) \cite{Xue3} proposed the two-stage method and bias-corrected empirical log-likelihood to obtain the estimators of the regression parameters and coefficient functions, established the asymptotic theory of the estimators and constructed the confidence regions of the regression parameters and pointwise confidence intervals for the coefficient functions.

Variable selection has received a lot of attention in statistical modeling and data analysis recently. %The development of variable selection has experienced rapidly since the 1970s.
Most of the variable selection procedures are based on penalized estimation using penalty functions, such as, $L_q$ penalty in Frank and Friedman (1993) \cite{Frank}, Lasso penalty in Tibshirani (1996) \cite{Tibshirani}, adaptive Lasso in Zou (2006) \cite{ZouH}, smoothly clipped absolute deviation(SCAD) penalty in Fan and Li (2001) \cite{FanJL}, among others. Fan and Li (2001) \cite{FanJL}
%has three excellent properties which are unbiasedness, sparsity and continuity.
indicated that SCAD penalty can not only select important variables consistently, but also produce the parameter estimators with the oracle property. An important issue of variable selection is how to choose the tuning parameter, some researchers studied different types of information criteria, for example, Akaike (1973) \cite{Akaike} proposed the Akaike information criterion (AIC), Schwarz (1978) \cite{Schwarz} investigated the Bayesian information criterion (BIC), and Foster and George (1994) \cite{Foster} introduced the risk inflation criterion (RIC). Pan (2001) \cite{PanW} developed the quasi-likelihood under the independence model criterion (QIC).

 In this paper, we consider the variable selection problem for model (\ref{1.1}).
  %{\color{blue} {we consider the variable selection problem for model (\ref{1.1})}}.
  %In particular, we will
 First, we use the B-spline functions to approximate the unknown coefficient functions in the model. Second, under the restriction of $\|\beta\|=1$, we employ the ``reparameterization" approach which has been applied by Yu and Ruppert (2002) \cite{Yu} and Wang et al. (2010) \cite{Wangj} to establish the penalized least-squares function. Third, we apply cross-validation method to select the tuning parameters and knots, and propose a stepwise iterative algorithm to compute the estimators. Under some regularization conditions, we show that this variable selection procedure is consistent and the estimators have the oracle property~i.e. sparsity and asymptotic normality, which means that the estimators of the parametric components have the same asymptotic distribution as that based on the correct submodel and the estimators of the nonparametric components achieve the optimal convergence rate. Compared with Feng and Xue (2015) \cite{Feng} and Zhao et al. (2019) \cite{ZhaoYX}, our method can select significant variables and estimate the regression parameters and  unknown coefficient functions simultaneously. This implies that our method can avoid the heavy computational burden.
 %{\color{red} {our method offers the following improvement which is that our method based on SCAD penalty This implies that our method can avoid the heavy computational burden.}}

The structure of the rest of this paper is presented as follows. In Section 2, we construct the penalized least-squares function using basis expansion, reparameterization and the SCAD penalty, and we also give some asymptotic properties about our variable selection approach, including the consistency of the variable selection and the
oracle property of the regularized estimators. In Section 3, we propose an iterative algorithm to find the penalized estimators based on local quadratic approximations and
how to select the tuning parameters. In Section 4, some Monte Carlo simulations and an
application using real data are carried out to evaluate the performance of the proposed method.
The proofs of the main results are provided in the Appendix.

\section{Variable selection via SCAD penalty}
\hspace{2em}Following the idea of Fan and Li (2001) \cite{FanJL}, we define the semiparametric penalized least-squares function as
\begin{align}
Q(\beta,\theta,\boldsymbol{g}(\cdot))=&\sum_{i=1}^n\{Y_i-\theta^{T}U_{i}-\boldsymbol{g}^\mathrm{T}(\beta^\mathrm{T}X_i)Z_i\}^2+n\sum_{l=1}^pp_{\lambda_{1l}}(|\beta_l|)\nonumber\\
&+n\sum_{h=1}^dp_{\lambda_{2h}}(|\theta_h|)
+n\sum_{k=1}^qp_{\lambda_{3k}}(\|g_k(\cdot)\|),\label{2.1}
\end{align}
where $\|g_k(\cdot)\|=\bigg(\int g_k^2(u)\mathrm{d}u\bigg)^{1/2}$, and $p_\lambda(\cdot)$ is the SCAD penalty function with a tuning parameter $\lambda$ which may
be chosen by a data-driven method. The first order derivative of ${p}_{\lambda}(\omega)$ is defined as
$$\dot{p}_{\lambda}(\omega)=\lambda\Big\{I(\omega\leqslant\lambda)+\frac{(a\lambda-\omega)_{+}}{(a-1)\lambda}\:I(\omega>\lambda)\Big\}$$
with $a>2,\omega>0$, and $p_\lambda(0)=0.$ We should state that the tuning parameters $\lambda_{1l}$, $\lambda_{2h}$ and $\lambda_{3k}$ are not necessarily the same for all $\beta_l$, $\theta_h$ and $g_k(\cdot)$.

Firstly, we can see that (\ref{2.1}) is unable to optimize since $\boldsymbol{g}(\cdot)$ is composed of unknown nonparametric functions. Similar to He et al. (2002) \cite{HeXM}, we consider basis function approximations for $\boldsymbol{g}(\cdot)$ in (\ref{2.1}). In particular, let $B(u)=\left(B_1(u),\ldots,B_L(u)\right)^\mathrm{T}$
be the B-spline basis functions with the order of $M+1$, where $L=K+M+1$, and $K$ is the number of interior knots. Thus, $g_k(u)$ can be approximated by
$$g_k(u)\approx B^\mathrm{T}(u)\gamma_k,\quad k=1,\ldots,q.$$
Substituting the estimators into (\ref{2.1}), we can obtain that
\begin{align}
Q(\beta,\theta,\gamma)=&\sum_{i=1}^{n}\{Y_{i}-\theta^{T}U_{i}-W_{i}^{\mathrm{T}}(\beta)\gamma\}^{2}+n\sum_{l=1}^{p}p_{\lambda_{1l}}(|\beta_{l}|)\nonumber\\
&+n\sum_{h=1}^{d}p_{\lambda_{2h}}(|\theta_{h}|)+n\sum_{k=1}^{q}p_{\lambda_{3k}}(\|\gamma_{k}\|_{H}),\label{2.2}
\end{align}
where $\gamma=(\gamma_{1}^{\mathrm{T}},\ldots,\gamma_{q}^{\mathrm{T}})^{\mathrm{T}},$ $W_{i}(\beta)=I_{q}\otimes B(\beta^{\mathrm{T}}X_{i})Z_{i},$ $\|\gamma_k\|_H=(\gamma_k^\mathrm{T}H\gamma_k)^{1/2},$ $H=\int B(u)B^\mathrm{T}(u)\mathrm{d}u.$

Secondly, the constraint $\|\beta\|=1$ means that $\beta$ is the boundary point on the unit sphere, each component of $\boldsymbol{g}(\beta^Tx)$  is not differentiable with respect to $\beta$. Therefore, we handle with ``reparameterization" approach. In particular, we write $\phi=(\beta_2,\ldots,\beta_p)$ which is one dimension lower than $\beta$, and define
$$\beta=\beta(\phi)=(\sqrt{1-\|\phi\|^2},~\phi^\mathrm{T})^\mathrm{T}.$$
%Then the true parameter $\phi_0$ must satisfy $\|\phi_0\|<1$, which is an inequality constraint.
Then $\beta(\phi)$ is infinitely differentiable with respect to $\phi$ under the constraint $\|\phi\|<1$. We can calculate the Jacobian matrix of $\beta$ with respect to $\phi$ by
\begin{align}\label{2.3}
J_\phi=\left(\begin{array}{c}-(1-\|\phi\|^2)^{-1/2}\phi^\mathrm{T}\\\mathbf{I}_{p-1}\end{array}\right),
\end{align}
where $\mathbf I_p$ is the $p\times p$ identity matrix, and the penalized least-squares function (\ref{2.2}) can be transformed to
\begin{align}
Q(\phi,\theta,\gamma)=&\sum_{i=1}^{n}\{Y_{i}-\theta^{T}U_{i}-W_{i}^{\mathrm{T}}(\phi)\gamma\}^{2}+n\sum_{l=1}^{p-1}p_{\lambda_{1l}}(|\phi_{l}|)\nonumber\\
&+n\sum_{h=1}^{d}p_{\lambda_{2h}}(|\theta_{h}|)+n\sum_{k=1}^{q}p_{\lambda_{3k}}(\|\gamma_{k}\|_{H}).\label{2.4}
\end{align}
Let $\hat{\phi}$, $\hat{\theta}$ and $\hat{\gamma}\:=\:(\hat{\gamma}_{1}^{\mathrm{T}},\ldots,\hat{\gamma}_{q}^{\mathrm{T}})^{\mathrm{T}}$ be the solution to the minimum of (\ref{2.4}). Then, the penalized least-squares estimator of $\beta$ is
$$\hat{\beta}=\left(\sqrt{1-\|\hat{\phi}\|^2},~\hat{\phi}^\mathrm{T}\right)^\mathrm{T},$$
and the estimator of $g_k(u)$ can be obtained by
$$\hat{g}_k(u)\approx B^\mathrm{T}(u)\hat{\gamma}_k.$$

Finally, we investigate the asymptotic properties of the penalized least-squares estimators. Let $\beta_0$, $\theta_0$, $\gamma_0$ and $\boldsymbol g_0(\cdot)$ be the true values of ${\beta}$, ${\theta}$, $\gamma$ and $\boldsymbol g(\cdot)$, respectively. We assume that $\beta_{l0}=0$, $l=s+1,\ldots,p$, and $\beta_{l0}$, $l=1,\ldots,s$ are all nonzero components of $\beta_0$; $\theta_{h0}=0$, $h=w+1,\ldots,d$, and $\theta_{h0}$, $l=1,\ldots,w$ are all nonzero components of $\theta_0$. Furthermore, we assume that $g_{k,0}(\cdot)=0$, $k=v+1,\ldots,q$, and $g_{k,0}(\cdot)$, $k=1,\ldots,v$ are all nonzero components of $\boldsymbol g_0(\cdot).$ The following theorem gives the consistency of the penalized least-squares estimators.
\begin{theorem}\label{th1}
~\rm{Suppose that the regularity conditions} \rm{(C1)}-\rm{(C7)} in the Appendix
hold and the number of knots $K= O_{p}( \emph{n}^{1/ ( 2\emph{r}+ 1) }).$ Then we have
\begin{flalign*}
&(\text{\rm i})\quad \| \hat{\beta } - \beta _0\| = O_{p}( n^{- r/ ( 2r+ 1) }+ a_n);&\\
&(\text{\rm ii})\quad \| \hat{\theta } - \theta _0\| = O_{p}( n^{- r/ ( 2r+ 1) }+ a_n);&\\
&(\text{\rm iii})\quad \| \hat{g} _k( \cdot ) - g_{k,0}( \cdot ) \| = O_{p}( n^{- r/ ( 2r+ 1) }+ a_n),~k= 1, \ldots , q,&
\end{flalign*}
where
$$a_{n}=\max_{l,h,k}\left\{|\dot{p}_{\lambda_{1l}}(|\beta_{l0}|)|,\:|\dot{p}_{\lambda_{2h}}(|\theta_{h0}|)|,\:|\dot{p}_{\lambda_{3k}}(\|\gamma_{k0}\|_{H})|\colon\:\beta_{l0}\neq0,\:\theta_{h0}\neq0,\:\gamma_{k0}\neq0\right\},$$
$r$ is defined in condition \rm{(C2)} in the Appendix, and $\dot{p} _{\lambda }( \cdot )$ denotes the first order derivative of $p_\lambda(\cdot).$
\end{theorem}

Furthermore, under some conditions, we show that the penalized least-squares estimators also have the sparsity property as follows.
\begin{theorem}\label{th2}
~\rm{Suppose that the regularity conditions \rm{(C1)}-\rm{(C7)} in the Appendix hold and the number of knots $K= O_{p}( \emph{n}^{1/ ( 2\emph{r}+ 1) }).$ Let
$$\lambda_{\mathrm{max}}=\max_{l,h,k}\{\lambda_{1l},\lambda_{2h},\lambda_{3k}\},\quad\lambda_{\mathrm{min}}=\min_{l,h,k}\{\lambda_{1l},\lambda_{2h},\lambda_{3k}\}.$$
If $\lambda_{\max}\to 0$  and $n^{r/(2r+1)}\lambda_{\min}\to\infty$ as $n\to\infty$, then, with probability tending to 1, the estimators $\hat{\beta}$,  $\hat{\theta}$ and $\boldsymbol{\hat{g}}(\cdot)$   satisfy}
\begin{flalign*}
&(\text{\rm i})\quad\hat{\beta}_l= 0,~l=s+1,\ldots,p;&\\
&(\text{\rm ii})\quad\hat{\theta}_h= 0,~h=w+1,\ldots,d;&\\
&(\text{\rm iii})\quad\hat{g}_k(\cdot)=0,~k=v+1,\ldots,q.&
\end{flalign*}
\end{theorem}

Next, we show that the estimators for nonzero coefficients in the parametric components have the same asymptotic distribution with the correct submodel. To the end, we need more notations to present the asymptotic property of the resulting estimators. Let
%~$\gamma^*=(\gamma_{1}^{\mathrm{T}},\ldots,\gamma_{v}^{\mathrm{T}})^{\mathrm{T}}$, $\gamma_0^*$ $\gamma^*$
$\beta^*=(\beta_1,\ldots,\beta_s)^\mathrm{T}$, $\theta^*=(\theta_1,\ldots,\theta_w)^\mathrm{T}$, $\gamma^*=(\gamma_{1}^{\mathrm{T}},\ldots,\gamma_{v}^{\mathrm{T}})^{\mathrm{T}}$ and $\boldsymbol{g}^*(u)=(g_1^\mathrm{T}(u),\ldots,g_v^\mathrm{T}(u))^\mathrm{T}$, and denote $\beta_0^*$, $\theta_0^*$, $\gamma_0^*$ and $\boldsymbol{g}_0^*(u)$ be the true values of $\beta^*$, $\theta^*$, $\gamma^*$ and $\boldsymbol{g}^*(u)$, respectively. Corresponding covariates are denoted by $X_i^*$, $U_i^*$ and $Z_i^*,$ $i=1,\ldots,n.$ In addition, we denote
${\Sigma=
\begin{pmatrix}
\Sigma_{11}&\Sigma_{12}\\
\Sigma_{21}&\Sigma_{22}
\end{pmatrix},}$
where
\begin{equation*}
\begin{aligned}
&\Sigma_{11}=E\left(V^{*} V^{*\mathrm{T}}\right)
-E\left\{C_{1}\left(\beta_{0}^{*\mathrm{T}} X^{*}\right) D^{-1}\left(\beta_{0}^{*\mathrm{T}} X^{*}\right) C_{1}^{\mathrm{T}}\left(\beta_{0}^{*\mathrm{T}} X^{*}\right)\right\},\\
&\Sigma_{12}=E\left(V^{*} U^{* \mathrm{T}}\right)
-E\left\{C_{1}\left(\beta_{0}^{*\mathrm{T}} X^{*}\right) D^{-1}\left(\beta_{0}^{*\mathrm{T}} X^{*}\right) C_{2}^{\mathrm{T}}\left(\beta_{0}^{*\mathrm{T}} X^{*}\right)\right\},\\
&\Sigma_{21}=E\left(U^{*} V^{*\mathrm{T}}\right)
-E\left\{C_{2}\left(\beta_{0}^{*\mathrm{T}} X^{*}\right) D^{-1}\left(\beta_{0}^{*\mathrm{T}} X^{*}\right) C_{1}^{\mathrm{T}}\left(\beta_{0}^{*\mathrm{T}} X^{*}\right)\right\},\\
&\Sigma_{22}=E\left(U^{*} U^{*\mathrm{T}}\right)
-E\left\{C_{2}\left(\beta_{0}^{*\mathrm{T}} X^{*}\right) D^{-1}\left(\beta_{0}^{*\mathrm{T}} X^{*}\right) C_{2}^{\mathrm{T}}\left(\beta_{0}^{*\mathrm{T}} X^{*}\right)\right\},
\end{aligned}
\end{equation*}
%${\Sigma=
%\begin{pmatrix}
%\Sigma_{11}&\Sigma_{12}\\
%\Sigma_{21}&\Sigma_{22}
%\end{pmatrix},}$
and $V^*$, $C_1(u)$, $C_2(u)$ and $D(u)$ are defined in condition (C7) in the Appendix. We assume that $\Sigma$ is an invertible matrix. The following result states the asymptotic normality of
$\begin{pmatrix}
\hat{\beta}^{*}\\
\hat{\theta}^{*}
\end{pmatrix}.$

\begin{theorem}\label{th3}
~\rm{Under the assumptions of Theorem \ref{th2}, we have}
$$\begin{aligned}
&\sqrt{n}
\begin{pmatrix}
\hat{\beta}^{*}-\beta_{0}^{*}\\
\hat{\theta}^{*}-\theta_{0}^{*}
\end{pmatrix}\stackrel{\mathcal{L}}{\longrightarrow}N(0,\sigma^{2}\widetilde{J}_{\phi_{0}^{*}}\Sigma^{-1}\widetilde{J}_{\phi_{0}^{*}}^{\mathrm{T}}),\end{aligned}$$
\rm{where `$\stackrel{\mathcal{L}}{\longrightarrow}$' represents the convergence in distribution and $\widetilde{J}_{\phi_{0}^{*}}=\begin{pmatrix}
  J_{\phi_{0}^{*}}& 0\\
 0 &I_w
\end{pmatrix}$.}
\end{theorem}

\begin{remark}
\rm{Theorem \ref{th1} indicates that our variable selection method is consistent and the estimators of nonparametric components achieve the optimal convergence rate as if the subset of true zero coefficients is already known (see Stone (1982) \cite{Stone}).}  Theorems \ref{th2} and \ref{th3} show that the penalized estimators have the oracle property.
\end{remark}
\section{Algorithm}
\hspace{2em}Because the SCAD-penalty function is single at the origin, it is not applicable to use gradient approach. In this section, we develop an iterative algorithm based on the local quadratic approximation of the penalty function $p_\lambda(\cdot)$ as in Fan and Li (2001) \cite{FanJL}.
Specifically, in a neighborhood of a given nonzero $\omega_0$, an approximation of the penalty function at value $\omega_0$ can be given by
$$p_{\lambda}(|\omega|)\approx p_{\lambda}(|\omega_{0}|)+\frac{1}{2}\:\frac{\dot{p}_{\lambda}(|\omega_{0}|)}{|\omega_{0}|}\:(\omega^{2}-\omega_{0}^{2}).$$
Hence, for the given initial value $\phi_l^{(0)}$ with $|\phi_l^{(0)}|>0$, $l=1,\ldots,p-1$, $\theta_h^{(0)}$ with $|\theta_h^{(0)}|>0$, $h=1,\ldots,d$, and $\gamma_k^{(0)}$ with $\| \gamma _k^{( 0) }\| _H> 0$, $k= 1, \ldots , q$, we have
$$p_{\lambda_{1l}}(|\phi_{l}|)\approx p_{\lambda_{1l}}(|\phi_{l}^{(0)}|)+\frac{1}{2}\:\frac{\dot{p}_{\lambda_{1l}}(|\phi_{l}^{(0)}|)}{|\phi_{l}^{(0)}|}\:(|\phi_{l}|^{2}-|\phi_{l}^{(0)}|^{2}),$$
$$p_{\lambda_{2h}}(|\theta_{h}|)\approx p_{\lambda_{2h}}(|\theta_{h}^{(0)}|)+\frac{1}{2}\:\frac{\dot{p}_{\lambda_{2h}}(|\theta_{h}^{(0)}|)}{|\theta_{h}^{(0)}|}\:(|\theta_{h}|^{2}-|\theta_{h}^{(0)}|^{2}),$$
$$p_{\lambda_{3k}}(\|\gamma_k\|_H)\approx p_{\lambda_{3k}}(\|\gamma_k^{(0)}\|_H)+\frac{1}{2}\frac{\dot{p}_{\lambda_{3k}}(\|\gamma_k^{(0)}\|_H)}{\|\gamma_k^{(0)}\|_H}(\|\gamma_k\|_H^2-\|\gamma_k^{(0)}\|_H^2).$$
Let
$\widetilde{W}_i^T(\phi)=(U_i^T,W^T_i(\phi))^T$, $\alpha=(\theta^T,\gamma^T)^T$,
\begin{align}
&\Sigma(\phi)=\mathrm{diag}\left\{\frac{\dot{p}_{\lambda_{11}}(|\phi_1^{}|)}{|\phi_1^{}|},\ldots,\frac{\dot{p}_{\lambda_{1,p-1}}(|\phi_{p-1}^{}|)}{|\phi_{p-1}^{}|}\right\},\nonumber\\
&\Sigma(\alpha)=\mathrm{diag}\left\{\frac{\dot{p}_{\lambda_{21}}(|\theta_1^{}|)}{|\theta_1^{}|},\ldots,\frac{\dot{p}_{\lambda_{2,d}}(|\theta_{d}^{}|)}{|\theta_{d}^{}|},\frac{\dot{p}_{\lambda_{31}}(\|\gamma_1^{}\|_H)}{\|\gamma_1^{}\|_H}H,\ldots,\frac{\dot{p}_{\lambda_{3q}}(\|\gamma_q^{}\|_H)}{\|\gamma_q^{}\|_H}H\right\}.\nonumber
\end{align}
Then, except for a constant term, (\ref{2.4}) becomes
\begin{equation}\label{3.1}
Q(\phi,\alpha)\approx\sum_{i=1}^n\{Y_i-\widetilde{W}_i^\mathrm{T}(\phi)\alpha\}^2+\frac{n}{2}\:\phi^\mathrm{T}\Sigma(\phi^{(0)})\phi+\frac{n}{2}\:\alpha^\mathrm{T}\Sigma(\alpha^{(0)})\alpha.
\end{equation}
With the aid of the local quadratic approximation, the Newton-Raphson algorithm can be applied to minimize the penalized least-squares function $Q(\phi,\alpha)$. In the following part, we propose a stepwise iterative Newton-Raphson algorithm to estimate the model parameters.

\noindent $ \rm \mathbf{ Step 1}$\quad Start with preliminary estimators  $\hat{\phi}^{(0)}$ and $\hat{\alpha}^{(0)}$. For example, the unpenalized estimators obtained by minimizing (\ref{2.4}) with $p_{\lambda}(\cdot)=0$ can be used.\\
\noindent $ \rm \mathbf{ Step 2}$\quad Calculate the penalized least-squares estimator by (\ref{3.1}) that
\begin{equation*}
\hat{\alpha}^{(m)}=\left(\sum_{i=1}^{n}\widetilde{W}_{i}(\hat{\phi}^{(0)})\widetilde{W}_{i}^{\mathrm{T}}(\hat{\phi}^{(0)})+\frac{n}{2}\Sigma(\hat{\alpha}^{(0)})\right)^{-1}\sum_{i=1}^{n}\widetilde{W}_{i}(\hat{\phi}^{(0)})Y_{i}.
\end{equation*}
\noindent $ \rm \mathbf{ Step 3}$\quad
Utilize the estimators $\hat{\alpha}^{(m)}$ obtained by Step 2, and minimize
\begin{equation*}
Q(\phi;\lambda_{1l})=\sum_{i=1}^n\{Y_i-\widetilde{W}_i^\mathrm{T}(\phi)\hat{\alpha}^{(m)}\}^2+\frac{n}{2}\:\phi^\mathrm{T}\Sigma(\hat{\phi}^{(0)})\phi.
\end{equation*}
We can get an estimator of $\phi$, say $\hat{\phi}^{(m)}$, and then obtain $\hat{\beta}^{(m)}$ via the transformation. During the iteration, once $|\hat{\phi}_l^{(m)}|$, $|\hat{\theta}_h^{(m)}|$, $\|\hat{\gamma}_k^{(m)}\|_H<\epsilon$, we set $\hat{\phi}_l^{(m)}=0$, $\hat{\theta}_h^{(m)}=0$, $\hat{\gamma}_k^{(m)}=0$, where $\epsilon>0$ is a small positive value. In our implementation, we select $\epsilon=10^{-2}$.\\
\noindent $ \rm \mathbf{ Step 4}$\quad Set $\hat{\phi^{(0)}}=\hat{\phi}^{(m)},\hat{\alpha}^{(0)}=\hat{\alpha}^{(m)}$, iterate Steps 2 and 3 until convergence, and denote the final estimators of $\phi$ and $\alpha$ as $\hat{\phi}$ and $\hat{\alpha}$. Then we can get the final estimators $\hat{\beta}$, $\hat{\theta}$ and $\hat{\gamma}$.

To implement this method, we should choose the number of interior knots $K$, and the tuning parameters $\lambda_{1l}$, $\lambda_{2h}$ and $\lambda_{3k}$. We use the similar cross-validation method to choose the tuning parameters as in Wang et al. (2008) \cite{LWang}. However, there are too many tuning parameters in our penalty functions, and the minimization problem for the cross-validation score over a higher-dimensional space is difficult. To overcome this difficulty, similar to Zhao and Xue (2009) \cite{ZHX}, we take the tuning parameters as
$$\lambda_{1l}=\frac{\lambda}{|\hat{\phi}_{l}^{u}|},\quad\lambda_{2h}=\frac{\lambda}{|\hat{\theta}_{h}^{u}|},\quad\lambda_{3k}=\frac{\lambda}{\|\hat{\gamma}_{k}^{u}\|_{H}},$$
where $\hat{\phi}_l^u$, $\hat{\theta}_h^u$ and $\hat{\gamma}_k^u$ are the unpenalized estimators of $\phi_l$, $\theta_h$ and $\gamma_k$, respectively. Thus we can select $\lambda$ and $K$ by minimizing the cross-validation score
\begin{equation*}
\mathrm{CV}(K,\lambda)=\sum_{i=1}^n\{Y_{i}-\hat{\theta}^{T}_{[i]}U_{i}-W_{i}^{\mathrm{T}}(\hat{\phi}_{[i]})\hat{\gamma}_{[i]}\}^2,
\end{equation*}
where  $\hat{\phi}_{[i]}$, $\hat{\theta}_{[i]}$ and $\hat{\gamma}_{[i]}$ are the solutions of (2.4) after deleting the $i$th subject.
\section{Numerical results}
\subsection{Simulation studies}
\hspace{2em}In this section, we conduct some Monte Carlo simulation studies to evaluate the performance of the proposed procedure. In each simulation, we generate random data from model (\ref{1.1}), where
$$\beta_0=(\beta_{1,0},\ldots,\beta_{10,0})^\mathrm{T},\quad\beta_{1,0}=\frac{1}{3},\quad\beta_{2,0}=\frac{2}{3},\quad\beta_{3,0}=\frac{2}{3},$$
$$\theta_0=(\theta_{1,0},\ldots,\theta_{10,0})^\mathrm{T},\quad\theta_{1,0}=2.0,\quad\theta_{2,0}=1.6\quad\theta_{3,0}=0.8,$$
$$\boldsymbol{g}_0(u)=\left(g_{1,0}(u),\ldots,g_{10,0}(u)\right)^\mathrm{T},\quad g_{1,0}(u)=2\cos(\pi u),\quad g_{2,0}(u)=1+3u^2.$$
While the remaining coefficients, corresponding to the irrelevant variables, are given by zeros. To perform this simulation, we take the covariates $Z_i$, $i=1,\ldots,n$, following the multivariate normal distribution with mean 0 and $Cov(Z_{i_k},Z_{i_l})=4\times0. 5^{| k- l| }$, $k,~l=1, \ldots , 10$. The covariates $U_i$, $i=1,\ldots,n$ follow the multivariate normal distribution with mean 0 and $Cov(U_{i_k},U_{i_l})=3\times0. 5^{| k- l| }$, $k,~l=1, \ldots , 10$. The covariates $X_i$, $i=1,\ldots,n$ are independent random vectors with each component uniformly distributed on (-0.75,~0.75). The model errors $\varepsilon_i$ will be considered under two independent $N( 0, 0.5^2)$ and $N( 0, 1.5^2)$ settings respectively.

We compare the performances of the variable selection procedures based on SCAD and LASSO penalty functions. In the following simulations, we use the cubic B-splines, and the sample size is set to $n=100, 150$ and $200$. The simulated results are reported in Tables 1-2 and Figures 1-2. The column labelled ``C" in Table 1 gives the average number of the true zero coefficients correctly set to zero, and the column labelled ``I" in Table 1 gives the average number of the true nonzero coefficients incorrectly set to zero.

($\rm{\Rmnum{1}}$)
In Tables 1-2, the mean, bias and standard deviation(SD) of the inner product $\hat{\beta}^T$$\beta_0$ are computed by 500 runs. The bias of the inner product $\hat{\beta}^T$$\beta_0$ is defined as
$$\text{\rm Bias}(\hat{\beta}^T\beta_0)=\hat{\beta}^T\beta_0-1.$$
The performance of estimator $\hat{\theta}$  will be assessed by using the generalized mean square error (GMSE), defined as
$$\mathrm{GMSE}(\hat{\theta})=(\hat{\theta}-\theta_0)^{\mathrm{T}}E(UU^{\mathrm{T}})(\hat{\theta}-\theta_0).$$

From the Tables 1-2, we can see these results:

(i) Under two different levels of model errors, the SCAD and LASSO variable selection methods both become closer to the oracle procedure as $n$ increases.  Specifically, the mean of the inner product $\hat{\beta}^T$$\beta_0$ is increasing, the bias and standard deviation(SD) of the inner product $\hat{\beta}^T$$\beta_0$ are decreasing as $n$ increases. Furthermore, the values in the column labeled ``C" of SCAD perform better than LASSO.

(ii) As the level of model errors decreases, the performance of SCAD becomes increasingly closer to that of the oracle procedure.

($\rm{\Rmnum{2}}$)
In Figures 1-2, the performances of estimator $\hat{{g}_1}(\cdot)$ and $\hat{{g}_2}(\cdot)$ will be exhibited by using the square root of
average square errors (RASE):
\begin{align}
\mathrm{RASE}_k=\left\{N^{-1}\sum_{j=1}^N[\hat{g}_k(u_j)-g_k(u_j)]^2\right\}^{1/2},\quad k=1,2,\nonumber
\end{align}
where $\{u_j,j=1,\ldots,N\}$ are the regular grid points at which the function $\hat{{g}_k}(u)$
is evaluated. In our simulation, $N=20$ is used. The performance of estimator $\hat{\boldsymbol{g}}(\cdot)$ will be assessed by $\rm{RASE}=\rm{RASE_1}+\rm{RASE_2}$.

From Figures 1-2, we can see that the estimated coefficient function curves of $g_1(u)$ and $g_2(u)$ by SCAD and LASSO are both close to the true coefficient function curves. As the level of model errors decreases, the estimated curves become slightly close to the true curves.
%\begin{table}[htbp]

%\caption{\label{tab:test}Variable selections based on adaptive tuning parameters(ATP) and constant tuning parameters (CTP)}

%\setlength{\tabcolsep}{1.4pt}
%\begin{adjustbox}{max width=\textwidth}
%\begin{tabular}{cccccccccc}
%\toprule
%&&&Estimator of $\beta$&&&&&Estimator of $\beta$&\\
%\cmidrule(r){3-7}\cmidrule(r){8-10}
%n   &Method   &Mean&Bias&SD&C&I&GMSE&C &I \\

%\midrule

%                  100&SCAD&0.9970 &-0.0078&0.1241&6.3460 &0.0180&0.01758&6.6560&0\\
%                     &LASSO&0.9971&-0.0079&0.1248&6.3440 &0.0160&0.01892&6.6340&0\\
%                  150&SCAD&0.9980&-0.0032&0.0791&6.6180&0&0.0098&6.7880&0\\
%                     &LASSO&0.9978&-0.0035&0.0796&6.5100&0 &0.0115&6.6920&0\\
%%                  200&SCAD&0.9995&-0.0006&0.0199&6.9320&0 &0.0066&6.9500&0\\
%                     &LASSO&0.9994&-0.0006&0.0214&6.8660&0 &0.0069&6.9160&0\\
%

%\bottomrule
%\end{tabular}
%\end{adjustbox}
%\end{table}
\begin{table}[H]

\caption{Simulation results for the estimator of $\beta$ and $\theta$($\varepsilon_i\sim N( 0, 0.5^2))$}

\setlength{\tabcolsep}{2.6pt}
%\begin{adjustbox}{max width=\textwidth}
\begin{tabular}{cccccccccc}
\toprule

&&\multicolumn{5}{c}{\textbf{\underline{~~~~~~~~~~~~~~~~~~Estimator of $\beta$~~~~~~~~~~~~~~}}}&\multicolumn{3}{c}{\textbf{\underline{~~~~~Estimator of $\theta$~~~~~}}}\\
\textbf{$n$}&\textbf{Method}&\textbf{Mean}&\textbf{Bias} &\textbf{SD} &\textbf{C}&\textbf{I}&\textbf{GMSE}&\textbf{C} &\textbf{I} \\

\midrule

                  100&SCAD&0.99522 &-0.00478&0.02335&6.45200 &0.00400&0.00938&6.71400&0\\
                     &LASSO&0.99578&-0.00422&0.02151&6.38200 &0.00400&0.01111&6.66800&0\\
                     &Oracle&0.99636&-0.00401&0.01987&      7&0      &0.00899&7      &0\\
                  150&SCAD&0.99826&-0.00174&0.00263&6.63800  &0.00200&0.00527&6.85400&0\\
                     &LASSO&0.99822&-0.00178&0.00232&6.54800 &0     &0.00633&6.80600&0\\
                     &Oracle&0.99835&-0.00170&0.00224&      7&0     &0.00512&7     &0\\
                  200&SCAD&0.99950&-0.00050 &0.00025 &6.98600&0     &0.00390&6.99000&0\\
                     &LASSO&0.99949&-0.00051&0.00025&6.95000 &0     &0.00476&6.96600&0\\
                     &Oracle&0.99958&-0.00049&0.00023&      7&0     &0.00387&7     &0\\

\bottomrule
\end{tabular}
%\end{adjustbox}
\end{table}

\begin{table}[htbp]

\caption{Simulation results for the estimator of $\beta$ and $\theta$($\varepsilon_i\sim N( 0, 1.5^2))$}

\setlength{\tabcolsep}{2.6pt}
%\begin{adjustbox}{max width=\textwidth}
\begin{tabular}{cccccccccc}
\toprule

&&\multicolumn{5}{c}{\textbf{\underline{~~~~~~~~~~~~~~~~~~Estimator of $\beta$~~~~~~~~~~~~~~}}}&\multicolumn{3}{c}{\textbf{\underline{~~~~~Estimator of $\theta$~~~~~}}}\\
\textbf{$n$}&\textbf{Method}&\textbf{Mean}&\textbf{Bias} &\textbf{SD} &\textbf{C}&\textbf{I}&\textbf{GMSE}&\textbf{C} &\textbf{I} \\

\midrule

                  100&SCAD&0.97324 &-0.02676&0.07658&6.22300  &0     &0.08017&6.18200&0\\
                     &LASSO&0.97176&-0.02824&0.08754&6.20000  &0     &0.08279&6.11600&0\\
                     &Oracle&0.97412&-0.02567&0.06387&      7&0     &0.07377&7      &0\\
                  150&SCAD&0.98886&-0.01114&0.01956&6.45100  &0     &0.04531&6.48600&0\\
                     &LASSO&0.98836&-0.01164&0.02056&6.44600 &0     &0.04775&6.44400&0\\
                     &Oracle&0.98987&-0.01022&0.01814&      7&0     &0.04456&7     &0\\
                  200&SCAD&0.99848&-0.00152 &0.00347 &6.81200&0     &0.03876&6.62600&0\\
                     &LASSO&0.99850&-0.00149&0.00339&6.80000 &0     &0.03913&6.59800&0\\
                     &Oracle&0.99911&-0.00146&0.00325&      7&0     &0.03747&7     &0\\

\bottomrule
\end{tabular}
%\end{adjustbox}
\end{table}

\begin{figure}[H]
	\centering
\captionsetup{justification=raggedright,singlelinecheck=false}
	\subfigure{
		%\begin{minipage}[t]{0.5\linewidth}
		\centering
		\includegraphics[height=2.1in,width=2.6in]{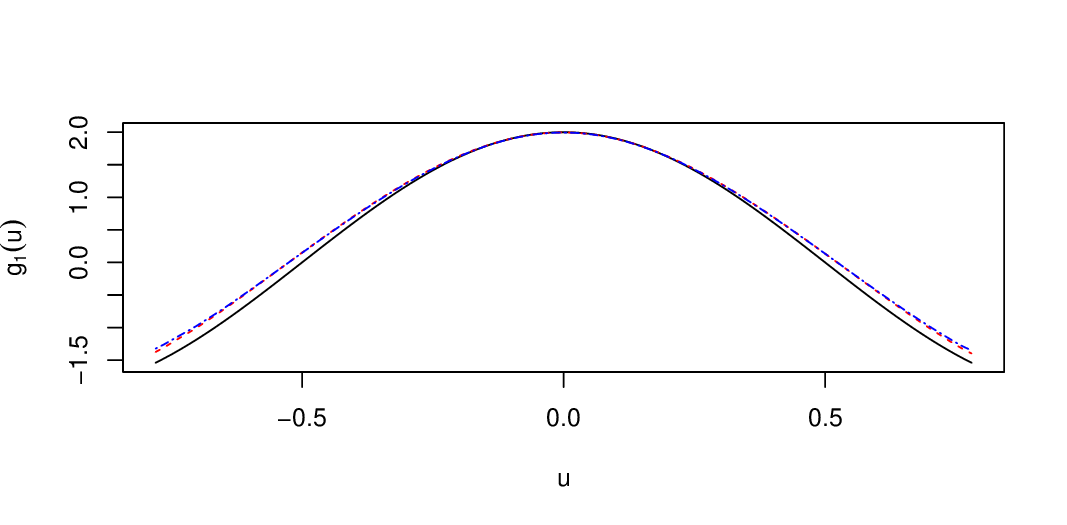}
		%\end{minipage}
	}
	\centering
	\subfigure{
		%\begin{minipage}[t]{0.5\linewidth}
		\centering
		\includegraphics[height=2.1in,width=2.6in]{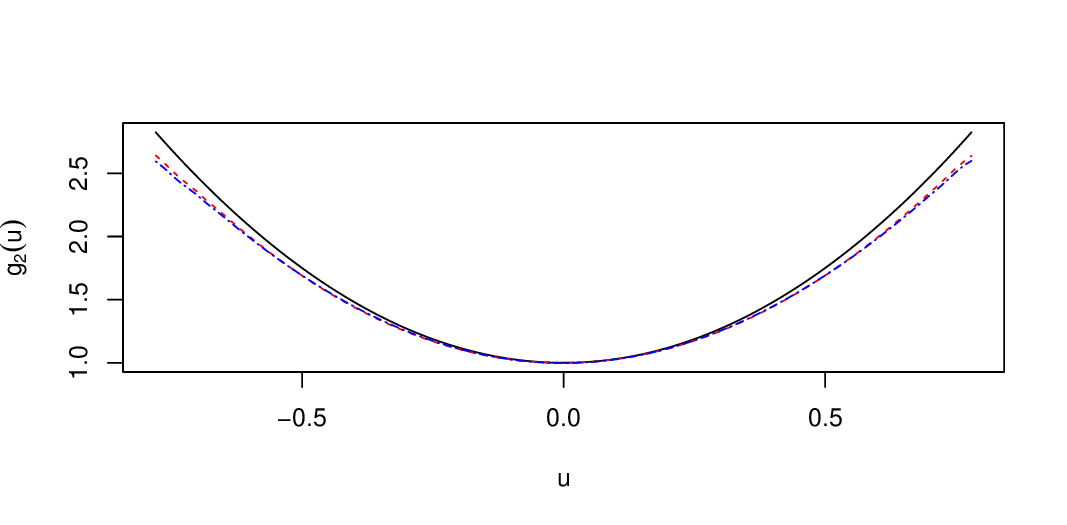}
	}
	\centering
	\caption{Simulation results when the $\varepsilon_i$ is independent $N( 0, 0.5^2)$ and $n=100$. The left panel shows the true curve(black solid curve), the SCAD estimated curve(red dashed curve) and the LASSO estimated curve(blue longdash curve) of $g_1(u)$. The right  panel shows the true curve(black solid curve), the SCAD estimated curve(red dashed curve) and the LASSO estimated curve(blue longdash curve) of $g_2(u)$.}
\end{figure}
\begin{figure}[H]
	\centering
\captionsetup{justification=raggedright,singlelinecheck=false}
	\subfigure{
		%\begin{minipage}[t]{0.5\linewidth}
		\centering
		\includegraphics[height=2.1in,width=2.6in]{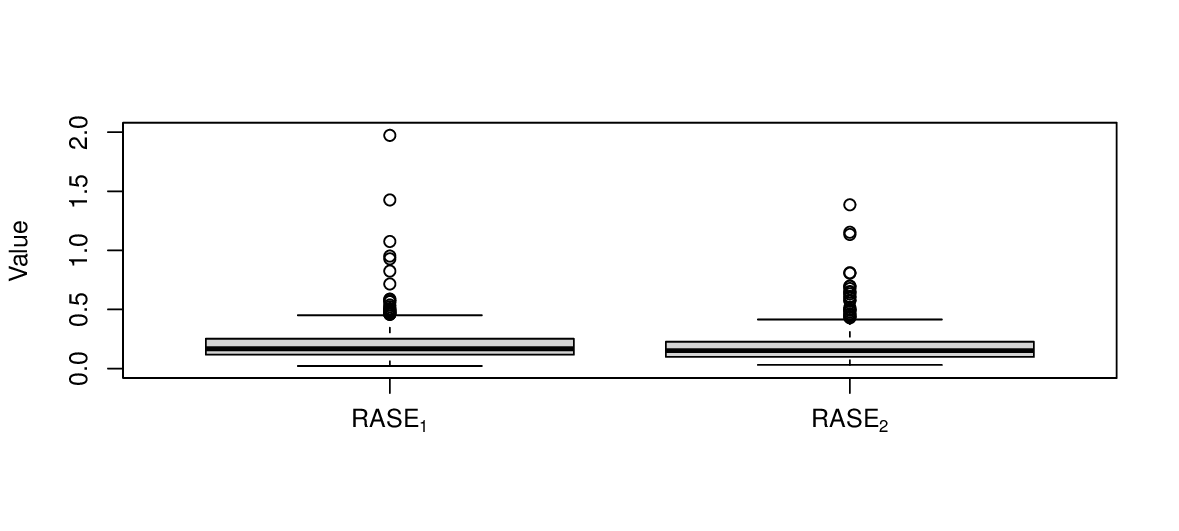}
		%\end{minipage}
	}
	\centering
	\subfigure{
		%\begin{minipage}[t]{0.5\linewidth}
		\centering
		\includegraphics[height=2.1in,width=2.6in]{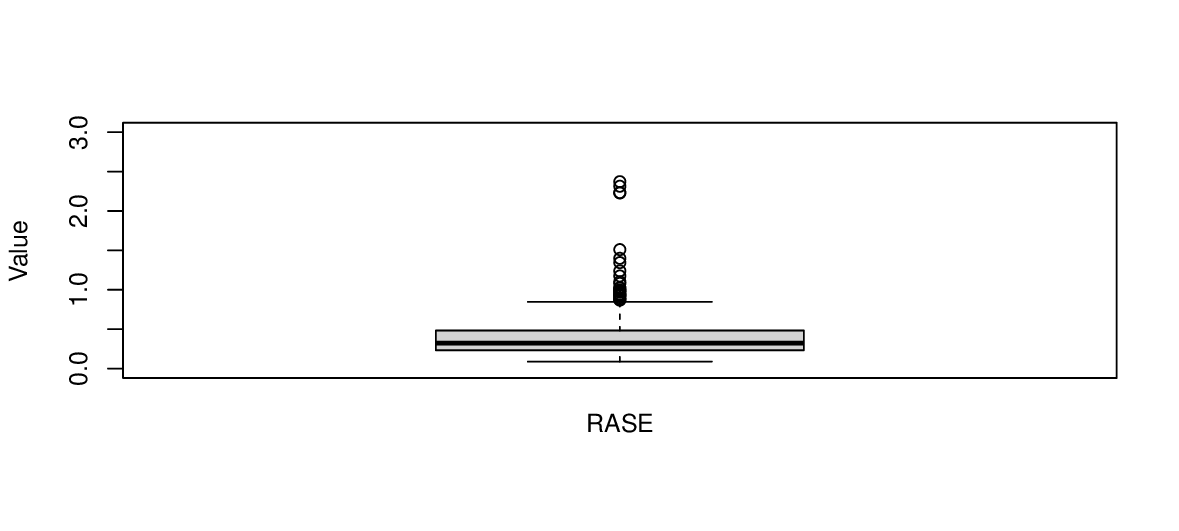}
	}
	\centering
	\caption{Simulation results when the $\varepsilon_i$ is independent $N( 0, 0.5^2)$ and $n=100$. The left panel is the boxplot the 500 $\rm{RASE_1}$ and $\rm{RASE_2}$ values of $g_1(u)$ and $g_2(u)$. The right panel is the boxplot the 500 RASE values of $\boldsymbol{g}(u)$.}
\end{figure}

\begin{figure}[H]
	\centering
\captionsetup{justification=raggedright,singlelinecheck=false}
	\subfigure{
		%\begin{minipage}[t]{0.5\linewidth}
		\centering
		\includegraphics[height=2.1in,width=2.6in]{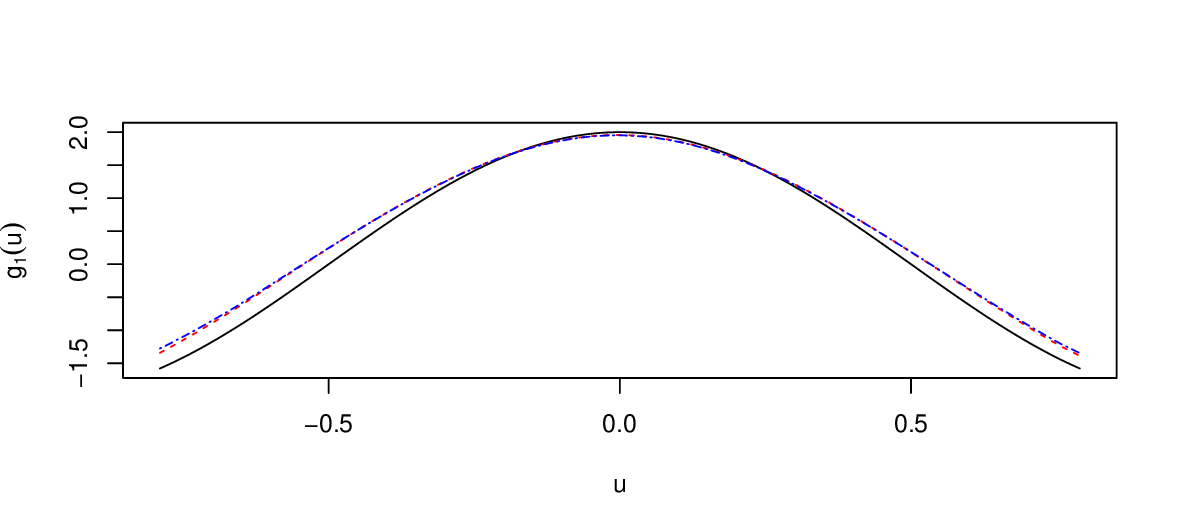}
		%\end{minipage}
	}
	\centering
	\subfigure{
		%\begin{minipage}[t]{0.5\linewidth}
		\centering
		\includegraphics[height=2.1in,width=2.6in]{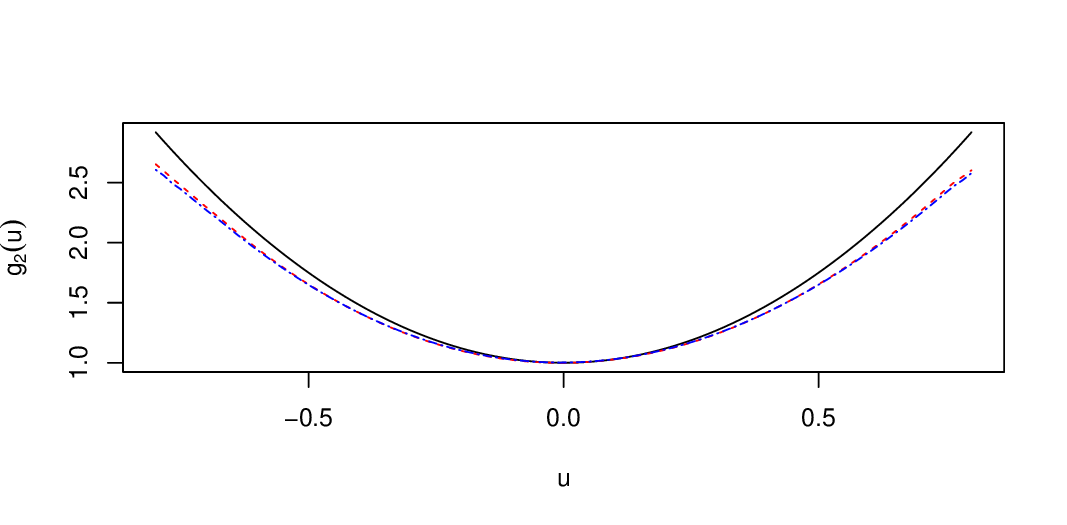}
	}
	\centering
	\caption{Simulation results when the $\varepsilon_i$ is independent $N( 0, 1.5^2)$ and $n=100$. The left panel shows the true curve(black solid curve), the SCAD estimated curve(red dashed curve) and the LASSO estimated curve(blue longdash curve) of $g_1(u)$. The right  panel shows the true curve(black solid curve), the SCAD estimated curve(red dashed curve) and the LASSO estimated curve(blue longdash curve) of $g_2(u)$.}
\end{figure}
\begin{figure}[H]
	\centering
\captionsetup{justification=raggedright,singlelinecheck=false}
	\subfigure{
		%\begin{minipage}[t]{0.5\linewidth}
		\centering
		\includegraphics[height=2.1in,width=2.6in]{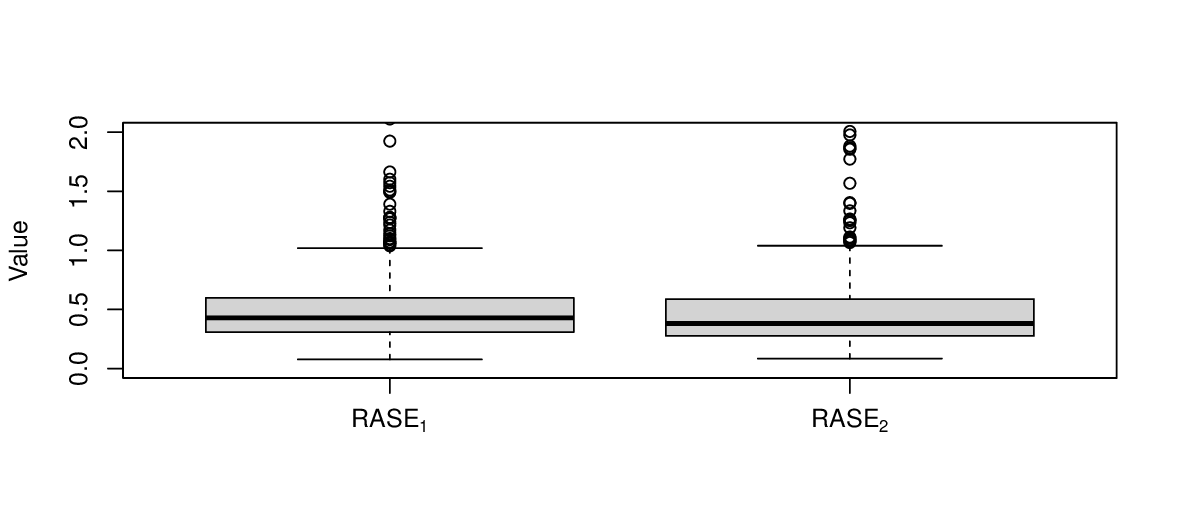}
		%\end{minipage}
	}
	\centering
	\subfigure{
		%\begin{minipage}[t]{0.5\linewidth}
		\centering
		\includegraphics[height=2.1in,width=2.6in]{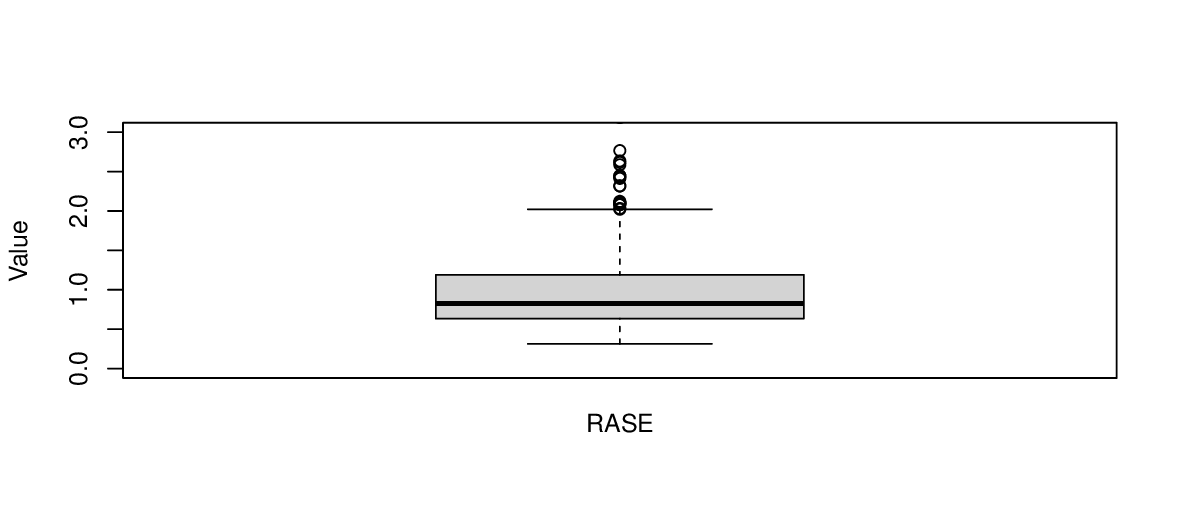}
	}
	\centering
	\caption{Simulation results when the $\varepsilon_i$ is independent $N( 0, 1.5^2)$ and $n=100$. The left panel is the boxplot the 500 $\rm{RASE_1}$ and $\rm{RASE_2}$ values of $g_1(u)$ and $g_2(u)$. The right panel is the boxplot the 500 RASE values of $\boldsymbol{g}(u)$.}
\end{figure}

%From Figure 2, we can see that the estimated coefficient function curves of $g_1(u)$ and $g_2(u)$ are both close to the true cofficient function curves by SCAD and LASSO. In particular, the estimated coefficient function curves based on SCAD is closer to the true cofficient function curves than LASSO. Figure 2 shows the RMSE value for the estimator of  cofficient function is very small. The results above indicate that the proposed method works well.

\subsection{Application to the body fat dataset}
\hspace{2em}We illustrate our proposed variable selection method by applying the body fat dataset, which is available from the website (http://lib.stat.cmu.edu/datasets/bodyfat). The response variable is the percentage of body fat which is determined by the underwater weighting technique (Siri, (1956) \cite{WE}). The thirteen covariates in this dataset are Age, Weight, Height, Neck, Chest, Abdomen, Hip, Thigh, Knee, Ankle, Biceps, Forearm and Wrist. The entire dateset contains 252 observations. Two observations are apparently type errors, four observations show inconsistency between the percentage of body fat and body density. These six observations can be removed from the dataset as outliers. In the following analysis, we employ the remain 246 observations.

In order to fit our model, we take logarithm of the percentage of body fat as the response variable $Y$, the covariates for linear part $U_1$=Age, $U_2$=Weight, the covariates dependent on varying coefficients $Z_1$=1, $Z_2$=Height, the index covariates $X_1$=Neck, $X_2$=Chest, $X_3$=Abdomen, $X_4$=Hip, $X_5$=Thigh, $X_6$=Knee, $X_7$=Ankle, $X_8$=Biceps, $X_9$=Forearm and $X_{10}$=Wrist. Then, we will establish the partially linear single-index variable-coefficient model (PLSIVM) as follows:
\begin{equation*}
Y=\theta_1U_{1}+\theta_2U_{2}+{g}(\beta^{T}X)Z_{1}+{g}(\beta^{T}X)Z_{2}+\varepsilon,
\end{equation*}
where $\beta^{T}X=\beta_1X_1+\cdots+\beta_{10}X_{10}$. For the PLSIVM, the restriction $\|\beta\|=1$ is required for identifiability.

From Table 3, we can obtain that:

(1) The proposed variable selection method chooses five covariates: Age($U_1$), Weight\\($U_2$), Neck($X_1$), Abdomen($X_3$), Wrist($X_{10}$). Among these significant covariates, Abdomen is the most important measurement for the prediction of the percentage of body fat, Wrist has more significant effects than Age, Weight and Neck, which is coincident with the results in  Peng and Huang (2011) \cite{Peng} and Lin and Huang (2013) \cite{LinH}.

(2) It is worth noting that Weight has nonzero constant on $Y$ which is similar to LM while is different from the existed results of SIM  in Peng and Huang (2011) \cite{Peng}, STM in Lin and Huang (2013) \cite{LinH} and SIVCM in Feng and Xue (2015) \cite{Feng}.

(3) The multiple $R^2$=0.70158 which is slightly larger than that of other four models. Therefore, our proposed variable selection method for PLSIVM is more appropriate to fit the body fat dataset.

\begin{table}[H]
\caption{Estimation results of body fat data}
\setlength\tabcolsep{12pt}
\begin{tabular}{ccccccccccc}
\toprule
                  Method&PLSIVCM& SIVCM    &SIM      &STM        &LM      \\
\midrule
                Age    & 0.00575& 0.00820 & 0.01490 & 0.08330   &0.04890 \\
                Weight & 0.00595& 0       & 0       & 0         &0.14570 \\
                Height & 0      & 0       & 0       & -0.10140  &-0.03950\\
                Neck   &-0.09926& -0.09900& -0.16910& -0.02800  &-0.14080\\
                Chest  &0       & 0       & 0       & 0         &-0.09430\\
                Abdomen&0.96980 & 0.96910 & 0.96060 & 0.97670   &0.76630 \\
                Hip    & 0      & 0       & 0       & 0         &-0.36380\\
                Thigh  & 0      & 0       & 0       & 0         &0.14610 \\
                Knee   & 0      & 0       & 0       & 0         &0       \\
                Ankle  & 0      & 0       & 0       & 0         &0       \\
                Biceps & 0      & 0       & 0       & 0         &0       \\
                Forearm& 0      & 0       & 0       & 0         &0.04130 \\
                Wrist  &-0.14629& -0.22600 & -0.22020 & -0.16770   &-0.11860\\
                       &        &         &         &           &       \\
                $R^2$  &0.70158 & 0.68160 & 0.67380 & 0.66650   &0.61480 \\
\bottomrule
\end{tabular}
\end{table}
\section*{Appendix.  Proof of theorems}
\hspace{2em}We begin the Appendix by listing some regularity conditions that are used in this paper. For convenience and simplicity, let $C$ denote a positive constant that may be different at each appearance throughout this paper.

(C1)\quad The density function $f(u)$ of $\beta_0^\mathrm{T}X,$ is bounded away from 0 on $\mathscr{U}=$ $\{u=\beta^{\mathrm{T}}x\colon x\in A\}$, where $A$ is the bounded support set of $X$. Furthermore, we assume that $f(u)$ satisfies the Lipschitz condition of order 1 on $\mathscr{U}.$

(C2)\quad The function $g_j(u)$, $j=1,\ldots,q$, have bounded and continuous derivatives up to order $r(\geqslant2)$ on $\mathscr{U}$, where $g_j(u)$ is the $j$th component of $\boldsymbol{g}(u).$

(C3)\quad$E(\|X\|^6)<\infty$, $E(\|U\|^6)<\infty$, $E(\|Z\|^6)<\infty$,  and $E(|\varepsilon|^6)<\infty.$

%(C4)\quad $\{ ( Y_i, X_i, U_i, Z_i) , 1\leqslant i\leqslant n\}$ is a strictly stationary and strongly mixing sequence with mixing coefficient $\alpha(n)=O(\rho^{n})$ for some $0<\rho<1$.

(C4)\quad Let $c_1,\ldots,c_K$ be the interior knots of $[a,b]$, where
$a=\inf\{u\colon u\in\mathscr{U}\},\quad b=\sup\{u\colon u\in\mathscr{U}\}.$
Furthermore, we let
$c_0=a,\quad c_{K+1}=b,\quad h_i=c_i-c_{i-1},\quad h=\max\{h_i\}.$
Then, there exists a constant $C_{0}$ such that
$$\frac{h}{\min\{h_{i}\}}<C_{0},\quad\max\{h_{i+1}-h_{i}\}=o(K^{-1}).$$

(C5)\quad $b_n\to 0$, as $n\to \infty$, where
$$b_{n}=\operatorname*{max}_{l,h,k}\left\{|\ddot{p}_{\lambda_{1l}}(|\phi_{l0}|)|, \:|\ddot{p}_{\lambda_{2h}}(|\theta_{h0}|), |\ddot{p}_{\lambda_{3k}}(\|\gamma_{k0}\|_{H})|\colon\phi_{l0}\neq0, \theta_{h0}\neq0,\:\gamma_{k0}\neq0\right\}.$$

(C6)\quad The penalty functions satisfy
\begin{align}
&\mathop{\lim\inf}_{n\to\infty}\mathop{\lim\inf}_{\phi_{l}\to0^{+}}\lambda_{1l}^{-1}|\dot{p}_{\lambda_{1l}}(|\phi_{l}|)|>0,\quad\mathop{\lim\inf}_{n\to\infty}\mathop{\lim\inf}_{\theta_{h}\to0^{+}}\lambda_{2h}^{-1}|\dot{p}_{\lambda_{2h}}(|\theta_{h}|)|>0,\nonumber\\ &\mathop{\lim\inf}_{n\to\infty}\mathop{\lim\inf}_{\|\gamma_{k}\|_{H}\to0}\lambda_{3k}^{-1}|\dot{p}_{\lambda_{3k}}(\|\gamma_{k}\|_{H})|>0,\nonumber
\end{align}
where $l= s, \ldots , p-1$, $h= w+1, \ldots , d$, $k= v+ 1, \ldots , q.$
%(C7)\quad
%$$\mathop{\lim\inf}_{n\to\infty}\mathop{\lim\inf}_{\phi_{l}\to0^{+}}\lambda_{1l}^{-1}|\dot{p}_{\lambda_{1l}}(|\phi_{l}|)|>0,\quad\mathop{\lim\inf}_{n\to\infty}\mathop{\lim\inf}_{\theta_{h}\to0^{+}}\lambda_{2h}^{-1}|\dot{p}_{\lambda_{2h}}(|\theta_{h}|)|>0,$$  %$$\mathop{\lim\inf}_{n\to\infty}\mathop{\lim\inf}_{\|\gamma_{k}\|_{H}\to0}\lambda_{3k}^{-1}|\dot{p}_{\lambda_{3k}}(\|\gamma_{k}\|_{H})|>0,$$
%where $l= s, \ldots , p-1$, $h= w, \ldots , d$, $k= v+ 1, \ldots , q.$

%(C7)\quad
%$\mathop{\lim\inf}_{n\to\infty}\mathop{\lim\inf}_{\phi_{l}\to0^{+}}\lambda_{1l}^{-1}|\dot{p}_{\lambda_{1l}}(|\phi_{l}|)|>0,~\mathop{\lim\inf}_{n\to\infty}\mathop{\lim\inf}_{\theta_{h}\to0^{+}}\lambda_{2h}^{-1}$\\
%$\cdot|\dot{p}_{\lambda_{2h}}(|\theta_{h}|)|>0,$ %~$\mathop{\lim\inf}_{n\to\infty}\mathop{\lim\inf}_{\|\gamma_{k}\|_{H}\to0}\lambda_{3k}^{-1}|\dot{p}_{\lambda_{3k}}(\|\gamma_{k}\|_{H})|>0,$\\
%where $l= s, \ldots , p-1$, $h= w, \ldots , d$, $k= v+ 1, \ldots , q.$

(C7)\quad The matrix
$D(u)=E\{Z^*Z^{*\text{T}}\mid\beta_0^{*\text{T}}X^*=u\}$
is positive definite, and each element of $D(u)$,
$C_{1}(u)=E\{V^{*}Z^{*\mathrm{T}}\mid\beta_{0}^{*\mathrm{T}}X^{*}=u\}$ and
$C_{2}(u)=E\{U^{*}Z^{*\mathrm{T}}\mid\beta_{0}^{*\mathrm{T}}X^{*}=u\}$
satisfies the Lipschitz condition of order 1 on $\mathscr{U}$, where
$$V^*=\dot{\boldsymbol{g}}_0^\mathrm{T}(\beta_0^\mathrm{*T}X^*)Z^*J_{\phi_0^*}^\mathrm{T}X^*.$$

\begin{remark}
\quad\rm Condition (C1) is used to bound the density function $f(u)$ of $\beta_0^\mathrm{T}X$ away from 0. Condition (C2) is the standard smoothness and continuity condition in the nonparametric estimation procedure. Condition (C3) is the necessary condition for the asymptotic properties of the estimators. Condition (C4) indicates that $c_0,\ldots,c_{K+1}$ is a $C_0$-quasi-uniform sequence of partitions of [0,1]. Condition (C5) and Condition (C6) are assumptions on the penalty functions. Condition (C7) ensures that the limiting variance for the estimator of $\beta_0$ and $\theta_0$ exists.
\end{remark}

To obtain the proofs of the theorems, the following lemmas are required.
\begin{lemma}\label{yl1}
\quad\rm If $g_k(u)$, $k= 1, \ldots , q$, satisfy condition (C2), then there exists $a$ constant $C> 0$ depending only on $M$ such that
\begin{equation*}
\sup\limits_{u\in\mathscr{U}}|g_k(u)-B^\mathrm{T}(u)\gamma_k|\leqslant CK^{-r}.
\end{equation*}
\end{lemma}
\noindent$\mathbf{Proof}$\quad This result is due to Corollary 6.21 in Schumaker (1981) \cite{Schumaker}, we omit the details here.$\hfill\square$\\
\begin{lemma}\label{yl2}
\quad\rm Suppose that conditions (C1)-(C3) and (C7) hold, and the number of knots $K= O( n^{- r/ ( 2r+ 1) }).$ Then we have
\begin{equation*}
\begin{aligned}
&\dfrac{1}{n}\sum_{i=1}^{n}V_{i}^{*}V_{i}^{*\text{T}}-\Psi_{n}^{\text{T}}\Phi_{n}^{-1}\Psi_{n}\xrightarrow{P}\Sigma_{11},\quad \dfrac{1}{n}\sum_{i=1}^{n}V_{i}^{*}U_{i}^{*\text{T}}-\Psi_{n}^{\text{T}}\Phi_{n}^{-1}G_{n}\xrightarrow{P}\Sigma_{12},\\
&\dfrac{1}{n}\sum_{i=1}^{n}U_{i}^{*}V_{i}^{*\text{T}}-G_{n}^{\text{T}}\Phi_{n}^{-1}\Psi_{n}\xrightarrow{P}\Sigma_{21},\quad
\dfrac{1}{n}\sum_{i=1}^{n}U_{i}^{*}U_{i}^{*\text{T}}-G_{n}^{\text{T}}\Phi_{n}^{-1}G_{n}\xrightarrow{P}\Sigma_{22},
\end{aligned}
\end{equation*}
and
\begin{equation*}
\begin{pmatrix}
 \frac{1}{n}\sum_{i=1}^{n}V_{i}^{*}V_{i}^{*\mathrm{T}}-\Psi_{n}^{\mathrm{T}}\Phi_{n}^{-1}\Psi_{n}
&\frac{1}{n}\sum_{i=1}^{n}V_{i}^{*}U_{i}^{*\mathrm{T}}-\Psi_{n}^{\mathrm{T}}\Phi_{n}^{-1}G_{n} \\  \frac{1}{n}\sum_{i=1}^{n}U_{i}^{*}V_{i}^{*\mathrm{T}}-G_{n}^{\mathrm{T}}\Phi_{n}^{-1}\Psi_{n}
&\frac{1}{n}\sum_{i=1}^{n}U_{i}^{*}U_{i}^{*\mathrm{T}}-G_{n}^{\mathrm{T}}\Phi_{n}^{-1}G_{n}
\end{pmatrix}\xrightarrow{P}\Sigma,
\end{equation*}
where ${\Sigma=
\begin{pmatrix}
\Sigma_{11}&\Sigma_{12}\\
\Sigma_{21}&\Sigma_{22}
\end{pmatrix},}$
$\Phi_{n}=\frac{1}{n}\sum_{i=1}^{n}W_{i}^{*}(\phi_{0}^{*})W_{i}^{*\mathrm{T}}(\phi_{0}^{*})$,\quad$\Psi_{n}=\frac{1}{n}\sum_{i=1}^{n}W_{i}^{*}(\phi_{0}^{*})V_{i}^{*\mathrm{T}}$,\quad $G_{n}=\frac{1}{n}\sum_{i=1}^{n}W_{i}^{*}(\phi_{0}^{*})U_{i}^{*\mathrm{T}}$, and `$\xrightarrow{P}$' means the convergence in probability.
\end{lemma}
\noindent$\mathbf{Proof}$\quad Let
$\mathbf{W}^{*}=(W_{1}^{*}(\phi_{0}^{*}),\ldots,W_{n}^{*}(\phi_{0}^{*}))^{\mathrm{T}}$, $\mathbf{V}^{*}=(V_{1}^{*},\ldots,V_{n}^{*})^{\mathrm{T}}$, $\mathbf{U}^{*}=(U_{1}^{*},\ldots,U_{n}^{*})^{\mathrm{T}}$,
$\mathbf{V}^{*}=\mathbf{V}^{*}-\Gamma_{n}+\Gamma_{n}=:\Delta_{n}+\Gamma_{n}$, $\mathbf{U}^{*}=\mathbf{U}^{*}-\Omega_{n}+\Omega_{n}=:\Xi_{n}+\Omega_{n},$
where
\begin{equation*}
\begin{aligned}
&V_{i}^{*}=\dot{\boldsymbol{g}}^{\mathrm{T}}(\beta_{0}^{*\mathrm{T}}X_{i}^{*})Z_{i}^{*}J_{\phi_{0}^{*}}^{\mathrm{T}}X_{i}^{*},\quad W_{i}^{*}(\phi_{0}^{*})=I_{q}\otimes B(\beta_{0}^{*\mathrm{T}}X_{i}^{*})Z_{i}^{*},\\
&\Gamma_{n}=(C_1(\beta_{0}^{*\mathrm{T}}X_{1}^{*})D^{-1}(\beta_{0}^{*\mathrm{T}}X_{1}^{*})Z_{1}^{*},\ldots,C_1(\beta_{0}^{*\mathrm{T}}X_{n}^{*})D^{-1}(\beta_{0}^{*\mathrm{T}}X_{n}^{*})Z_{n}^{*})^{\mathrm{T}},\\
&\Omega_{n}=(C_2(\beta_{0}^{*\mathrm{T}}X_{1}^{*})D^{-1}(\beta_{0}^{*\mathrm{T}}X_{1}^{*})Z_{1}^{*},\ldots,C_2(\beta_{0}^{*\mathrm{T}}X_{n}^{*})D^{-1}(\beta_{0}^{*\mathrm{T}}X_{n}^{*})Z_{n}^{*})^{\mathrm{T}}.
\end{aligned}
\end{equation*}

Then a simple calculation yields
\begin{align}
&\frac{1}{n}\sum_{i=1}^{n}V_{i}^{*}V_{i}^{*\mathrm{T}}-\Psi_{n}^{\mathrm{T}}\Phi_{n}^{-1}\Psi_{n}\nonumber\\
=~&n^{-1}\mathbf{V}^{*\mathrm{T}}(I-P^{\mathrm{T}})(I-P)\mathbf{V}^{*}\nonumber\\
=~&n^{-1}(\Delta_{n}+\Gamma_{n})^{\mathrm{T}}(I-P^{\mathrm{T}})(I-P)(\Delta_{n}+\Gamma_{n})\nonumber\\
=~&n^{-1}\{\Delta_{n}^{\mathrm{T}}\Delta_{n}
+\Delta_{n}^{\mathrm{T}}(I-P^{\mathrm{T}})(I-P)\Gamma_{n}
+\Gamma_{n}^{\mathrm{T}}(I-P^{\mathrm{T}})(I-P)\Delta_{n}\nonumber\\
&+\Gamma_{n}^{\mathrm{T}}(I-P^{\mathrm{T}})(I-P)\Gamma_{n}
-\Delta_{n}^{\mathrm{T}}P^{\mathrm{T}}P\Delta_{n}\},\label{lmgs1}
\end{align}
\begin{align}
&\frac{1}{n}\sum_{i=1}^{n}V_{i}^{*}U_{i}^{*\mathrm{T}}-\Psi_{n}^{\mathrm{T}}\Phi_{n}^{-1}G_{n}\nonumber\\
=~&n^{-1}\mathbf{V}^{*\mathrm{T}}(I-P^{\mathrm{T}})(I-P)\mathbf{U}^{*}\nonumber\\
=~&n^{-1}(\Delta_{n}+\Gamma_{n})^{\mathrm{T}}(I-P^{\mathrm{T}})(I-P)(\Xi_{n}+\Omega_{n})\nonumber\\
=~&n^{-1}\{\Delta_{n}^{\mathrm{T}}\Xi_n
+\Delta_{n}^{\mathrm{T}}(I-P^{\mathrm{T}})(I-P)\Omega_{n}
+\Gamma_{n}^{\mathrm{T}}(I-P^{\mathrm{T}})(I-P)\Xi_{n}\nonumber\\
&+\Gamma_{n}^{\mathrm{T}}(I-P^{\mathrm{T}})(I-P)\Omega_{n}
-\Delta_{n}^{\mathrm{T}}P^{\mathrm{T}}P\Xi_{n}\},\label{lmgs2}
\end{align}
\begin{align}
&\frac{1}{n}\sum_{i=1}^{n}U_{i}^{*}V_{i}^{*\mathrm{T}}-G_{n}^{\mathrm{T}}\Phi_{n}^{-1}\Psi_{n}\nonumber\\
=~&n^{-1}\mathbf{U}^{*\mathrm{T}}(I-P^{\mathrm{T}})(I-P)\mathbf{V}^{*}\nonumber\\
=~&n^{-1}(\Xi_{n}+\Omega_{n})^{\mathrm{T}}(I-P^{\mathrm{T}})(I-P)(\Delta_{n}+\Gamma_{n})\nonumber\\
=~&n^{-1}\{\Xi_n^{\mathrm{T}}\Delta_{n}
+\Xi_{n}^{\mathrm{T}}(I-P^{\mathrm{T}})(I-P)\Gamma_{n}
+\Omega_{n}^{\mathrm{T}}(I-P^{\mathrm{T}})(I-P)\Delta_{n}\nonumber\\
&+\Omega_{n}^{\mathrm{T}}(I-P^{\mathrm{T}})(I-P)\Gamma_{n}
-\Xi_{n}^{\mathrm{T}}P^{\mathrm{T}}P\Delta_{n}\},\label{lmgs3}
\end{align}
and
\begin{align}
&\frac{1}{n}\sum_{i=1}^{n}U_{i}^{*}U_{i}^{*\mathrm{T}}-G_{n}^{\mathrm{T}}\Phi_{n}^{-1}G_{n}\nonumber\\
=~&n^{-1}\mathbf{U}^{*\mathrm{T}}(I-P^{\mathrm{T}})(I-P)\mathbf{U}^{*}\nonumber\\
=~&n^{-1}(\Xi_{n}+\Omega_{n})^{\mathrm{T}}(I-P^{\mathrm{T}})(I-P)(\Xi_{n}+\Omega_{n})\nonumber\\
=~&n^{-1}\{\Xi_{n}^{\mathrm{T}}\Xi_{n}
+\Xi_{n}^{\mathrm{T}}(I-P^{\mathrm{T}})(I-P)\Omega_{n}
+\Omega_{n}^{\mathrm{T}}(I-P^{\mathrm{T}})(I-P)\Xi_{n}\nonumber\\
&+\Omega_{n}^{\mathrm{T}}(I-P^{\mathrm{T}})(I-P)\Omega_{n}
-\Xi_{n}^{\mathrm{T}}P^{\mathrm{T}}P\Xi_{n}\},\label{lmgs4}
\end{align}
where
$P=\mathbf{W}^*(\mathbf{W}^{*\mathrm{T}}\mathbf{W}^*)^{-1}\mathbf{W}^{*\mathrm{T}}.$

By  the proof of Lemma 1 in Zhao and Xue (2010) \cite{ZHX}, we easily obtain that all but the first term of (\ref{lmgs1}), (\ref{lmgs2}), (\ref{lmgs3}) and (\ref{lmgs4}) are $o_p(1)$. Then, by the law of large numbers, we can derive that the first term converges to $\Sigma_{11}$, $\Sigma_{12}$, $\Sigma_{21}$ and $\Sigma_{22}$ respectively. Hence, we have \begin{equation*}
\begin{pmatrix}
 \frac{1}{n}\sum_{i=1}^{n}V_{i}^{*}V_{i}^{*\mathrm{T}}-\Psi_{n}^{\mathrm{T}}\Phi_{n}^{-1}\Psi_{n}
&\frac{1}{n}\sum_{i=1}^{n}V_{i}^{*}U_{i}^{*\mathrm{T}}-\Psi_{n}^{\mathrm{T}}\Phi_{n}^{-1}G_{n} \\  \frac{1}{n}\sum_{i=1}^{n}U_{i}^{*}V_{i}^{*\mathrm{T}}-G_{n}^{\mathrm{T}}\Phi_{n}^{-1}\Psi_{n}
&\frac{1}{n}\sum_{i=1}^{n}U_{i}^{*}U_{i}^{*\mathrm{T}}-G_{n}^{\mathrm{T}}\Phi_{n}^{-1}G_{n}
\end{pmatrix}\xrightarrow{P}\Sigma.
\end{equation*}$\hfill\square$\\
$\mathbf{Proof~of~Theorem~1}$\quad Let $\delta=n^{-r/(2r+1)}+a_{n}$,\quad$\phi=\phi_{0}+\delta\tau_{1}$,\quad$\theta=\theta_{0}+\delta\tau_{2}$,\quad$\gamma=\gamma_{0}+\delta\tau_{3}$,
\quad$\tau=(\tau_{1}^{\mathrm{T}}$,$\tau_{2}^{\mathrm{T}}$,$\tau_{3}^{\mathrm{T}})^{\mathrm{T}}.$

(i) We will show that, for any given $\epsilon>0$, there exists a large constant $C$ such that

\begin{equation}
\begin{aligned}
P\Big\{\inf\limits_{\|\tau\|=C}Q(\phi,\theta,\gamma)>Q(\phi_0,\theta_0,\gamma_0)\Big\}\geqslant1-\epsilon,\label{mubiao1}
\end{aligned}
\end{equation}
where $\phi_0$, $\theta_0$ and $\gamma_0$ are true value of $\phi$, $\theta$ and $\gamma$.

Let
$L_n(\tau)=K^{-1}\{Q(\phi,\theta,\gamma)-Q(\phi_0,\theta_0,\gamma_0)\}.$  Thus, by the Taylor expansion and a direct calculation, we have
\begin{align}L_{n}(\tau)
=~&K^{-1}\{Q(\phi_{0}+\delta\tau_{1},\theta_{0}+\delta\tau_{2},\gamma_{0}+\delta\tau_{3})-Q(\phi_{0},\theta_{0},\gamma_{0})\}\nonumber\\
\geqslant~&\frac{-2\delta}{K}\sum_{i=1}^{n}(Y_{i}-\theta_0^{\mathrm{T}}U_{i}-W_{i}(\phi_0)^{\mathrm{T}}\gamma_0)(\dot{W}_{i}^{\mathrm{T}}(\phi_{0})\gamma_{0}\tau_{1}^{\mathrm{T}}J_{\phi_{0}}^{\mathrm{T}}X_{i}+U^{\mathrm{T}}\tau_{2}+W_{i}^{\mathrm{T}}(\phi_{0})\tau_{3})\nonumber\\
&+\frac{\delta^{2}}{K}\sum_{i=1}^{n}(\dot{W}_{i}^{\mathrm{T}}(\phi_{0})\gamma_{0}\tau_{1}^{\mathrm{T}}J_{\phi_{0}}^{\mathrm{T}}X_{i}+U^{\mathrm{T}}\tau_{2}+W_{i}^{\mathrm{T}}(\phi_{0})\tau_{3})^{2}\nonumber\\
&+\frac{n}{K}\sum_{l=1}^{s-1}[p_{\lambda_{1l}}(|\phi_{l}|)-p_{\lambda_{1l}}(|\phi_{l0}|)]\nonumber\\
&+\frac{n}{K}\sum_{h=1}^{w}[p_{\lambda_{2h}}(|\theta_{h}|)-p_{\lambda_{2h}}(|\theta_{h0}|)]\nonumber\\
&+\frac{n}{K}\sum_{k=1}^{v}[p_{\lambda_{3k}}(\|\gamma_{k}\|_{H})-p_{\lambda_{3k}}(\|\gamma_{k0}\|_{H})]+o_{p}(1),\nonumber
\end{align}
Let $\dot{W}_i(\phi_0)=I_q\otimes\dot{B}(\beta_0^\mathrm{T}X_i)Z_i,$\quad$R(u)=\left(R_{1}(u),\ldots,R_{q}(u)\right)^{\mathrm{T}}$,\quad$R_{k}(u)=g_{k0}(u)-B^{\mathrm{T}}(u)\gamma_{k0},$ $k=1,\ldots,q.$
Since $Y_{i}=\theta_0^{T}U_{i}+\boldsymbol{g}_0^{T}(\beta_0^{T}X_{i})Z_{i}+\varepsilon_{i}$, we can obtain that
$$\begin{aligned}
&Y_{i}-\theta_0^{\mathrm{T}}U_{i}-W_{i}(\phi_0)^{\mathrm{T}}\gamma_0\\ =~&\theta_0^{T}U_{i}+\boldsymbol{g}_0^{T}(\beta_0^{T}X_{i})Z_{i}+\varepsilon_{i}-\theta_0^{\mathrm{T}}U_{i}-W_{i}(\phi_0)^{\mathrm{T}}\gamma_0\\
=~&\varepsilon_{i}+\boldsymbol{g}_0^{T}(\beta_0^{T}X_{i})Z_{i}-W_{i}(\phi_0)^{\mathrm{T}}\gamma_0\\
=~&\varepsilon_{i}+R^{\mathrm{T}}(\beta_{0}^{\mathrm{T}}X_{i})Z_{i}.
\end{aligned}$$
Then we can write
\begin{align}
L_n(\tau)\geqslant~&\frac{-2\delta}{K}\sum_{i=1}^{n}(\varepsilon_{i}+R^{\mathrm{T}}(\beta_{0}^{\mathrm{T}}X_{i})Z_{i})(\dot{W}_{i}^{\mathrm{T}}(\phi_{0})\gamma_{0}\tau_{1}^{\mathrm{T}}J_{\phi_{0}}^{\mathrm{T}}X_{i}+U^{\mathrm{T}}\tau_{2}+W_{i}^{\mathrm{T}}(\phi_{0})\tau_{3})\nonumber\\
&+\frac{\delta^{2}}{K}\sum_{i=1}^{n}(\dot{W}_{i}^{\mathrm{T}}(\phi_{0})\gamma_{0}\tau_{1}^{\mathrm{T}}J_{\phi_{0}}^{\mathrm{T}}X_{i}+U^{\mathrm{T}}\tau_{2}+W_{i}^{\mathrm{T}}(\phi_{0})\tau_{3})^{2}+o_{p}(1)\nonumber\\
&+\frac{n}{K}\sum_{l=1}^{s-1}[p_{\lambda_{1l}}(|\phi_{l}|)-p_{\lambda_{1l}}(|\phi_{l0}|)]\nonumber\\
&+\frac{n}{K}\sum_{h=1}^{w}[p_{\lambda_{2h}}(|\theta_{h}|)-p_{\lambda_{2h}}(|\theta_{h0}|)]\nonumber\\
&+\frac{n}{K}\sum_{k=1}^{v}[p_{\lambda_{3k}}(\|\gamma_{k}\|_{H})-p_{\lambda_{3k}}(\|\gamma_{k0}\|_{H})],\nonumber\\
=:~&S_{1}+S_{2}+S_{3}+S_{4}+S_{5}+o_{p}(1).\label{jgS}
\end{align}
%ÁíÆðÒ»¶Î

Next, we analyse each term of (\ref{jgS}). For $S_1$, by conditions (C1), (C2), (C4), and Lemma \ref{yl1}, we can derive that
$$\|R_k(u)\|=O(K^{-r})$$
and
\begin{equation}
\begin{aligned}
|\dot{g}_{k}(\beta_{0}^{\mathrm{T}}X_{i})-\dot{B}^{\mathrm{T}}(\beta_{0}^{\mathrm{T}}X_{i})\gamma_{k0}|\leqslant CK^{-r+1}.\label{tiaojian}
\end{aligned}
\end{equation}
Then, a simple calculation yields

\begin{align}
&\sum_{i=1}^{n}R^{\mathrm{T}}(\beta_{0}^{\mathrm{T}}X_{i})Z_{i}\{\dot{W}_{i}^{\mathrm{T}}(\phi_{0})\gamma_{0}\tau_{1}^{\mathrm{T}}J_{\phi_{0}}^{\mathrm{T}}X_{i}+U^{\mathrm{T}}\tau_{2}+W_{i}^{\mathrm{T}}(\phi_{0})\tau_{3}\}\nonumber\\
=~&\sum_{i=1}^{n}R^{\mathrm{T}}(\beta_{0}^{\mathrm{T}}X_{i})Z_{i}\{\dot{\boldsymbol{g}}^{\mathrm{T}}(\beta_{0}^{\mathrm{T}}X_{i})Z_{i}\tau_{1}^{\mathrm{T}}J_{\phi_{0}}^{\mathrm{T}}X_{i}\nonumber\\
&+(\dot{W}_{i}^{\mathrm{T}}(\phi_{0})\gamma_{0}-\dot{\boldsymbol{g}}^{\mathrm{T}}(\beta_{0}^{\mathrm{T}}X_{i})Z_{i})\tau_{1}^{\mathrm{T}}J_{\phi_{0}}^{\mathrm{T}}X_{i}+U^{\mathrm{T}}\tau_{2}+W_{i}^{\mathrm{T}}(\phi_{0})\tau_{3}\}\nonumber\\
=~&O_p(nK^{-r}\|\tau\|).\label{4-4}
\end{align}
Note that $\varepsilon_i$ is independent of $(X_i,U_i,Z_i)$. We can prove that
$$\frac{1}{\sqrt{n}}\sum_{i=1}^{n}\varepsilon_{i}\{\dot{W}_{i}^{\mathrm{T}}(\phi_{0})\gamma_{0}\tau_{1}^{\mathrm{T}}J_{\phi_{0}}^{\mathrm{T}}X_{i}+U^{\mathrm{T}}\tau_{2}+W_{i}^{\mathrm{T}}(\phi_{0})\tau_{3}\}=O_{p}(\|\tau\|).$$
Combing this with (\ref{4-4}), we can get that
$$S_{1}=O_{p}(\sqrt{n}\:K^{-1}\delta)\|\tau\|+O_{p}(nK^{-1-r}\delta)\|\tau\|=O_{p}(1+n^{r/(2r+1)}a_{n})\|\tau\|.$$
%ÁíÆðÒ»¶Î

Similarly, we can prove that
$$S_2=O_p(nK^{-1}\delta^2)\|\tau\|^2=O_p(1+2n^{r/(2r+1)}a_n)\|\tau\|^2.$$
Hence, by choosing a sufficiently large $C$, $S_2$ dominates $S_1$ uniformly in $\|\tau\|=C.$
%ÁíÆðÒ»¶Î

By the fact that $p_\lambda(0)=0$ and the Taylor expansion, we get
$$\begin{aligned}S_{3}&\leqslant\frac{n}{K}\sum_{l=1}^{s-1}[\delta\dot{p}_{\lambda_{1l}}(|\phi_{l0}|)\mathrm{sgn}(\phi_{l0})|\tau_{1l}|+\delta^{2}\ddot{p}_{\lambda_{1l}}(|\phi_{l0}|)|\tau_{1l}|^{2}\{1+o(1)\}]\\&\leqslant\sqrt{s-1}\:K^{-1}n\delta a_{n}\|\tau\|+nK^{-1}\delta^{2}b_{n}\|\tau\|^{2},\end{aligned}$$
$$\begin{aligned}S_{4}&\leqslant\frac{n}{K}\sum_{l=1}^{w}[\delta\dot{p}_{\lambda_{2h}}(|\theta_{h0}|)\mathrm{sgn}(\theta_{h0})|\tau_{2h}|+\delta^{2}\ddot{p}_{\lambda_{2h}}(|\theta_{h0}|)|\tau_{2h}|^{2}\{1+o(1)\}]\\&\leqslant\sqrt{w}\:K^{-1}n\delta a_{n}\|\tau\|+nK^{-1}\delta^{2}b_{n}\|\tau\|^{2},\end{aligned}$$
$$\begin{aligned}S_{5}&\leqslant\frac{n}{K}\sum_{l=1}^{v}[\delta\dot{p}_{\lambda_{3k}}(|\gamma_{k0}|)\mathrm{sgn}(\gamma_{k0})|\tau_{3k}|+\delta^{2}\ddot{p}_{\lambda_{3k}}(|\gamma_{k0}|)|\tau_{3k}|^{2}\{1+o(1)\}]\\&\leqslant\sqrt{v}\:K^{-1}n\delta a_{n}\|\tau\|+nK^{-1}\delta^{2}b_{n}\|\tau\|^{2}.\end{aligned}$$
Then, it is easy to show that $S_3$, $S_4$ and $S_5$ are dominated by $S_2$ uniformly in $\|\tau\|=C.$  Hence, for a sufficiently large $C$, the probability of $L_n(\tau)=K^{-1}\{Q(\phi,\theta,\gamma)-Q(\phi_0,\theta_0,\gamma_0)\}>0$ is at least $1-\epsilon$, then, we can reach to (\ref{mubiao1}). Therefore, there exist local minimizers $\hat{\phi}$, $\hat{\theta}$ and $\hat{\gamma}$ such that
$$\|\hat\phi-\phi_0\|=O_p(\delta),\quad\|\hat\theta-\theta_0\|=O_p(\delta),\quad\|\hat\gamma-\gamma_0\|=O_p(\delta).$$
$\begin{array}{l}\text{A direct calculation can lead to }\|\hat{\beta}-\beta_{0}\|=O_{p}(\delta),\mathrm{~which~completes~the~proofs}\\\text{of (i).}\end{array}$

 (ii) It can be proved directly from the proof of (i).

(iii) Note that
$$\begin{aligned}\|\hat{g}_{k}(u)-g_{k0}(u)\|^{2}&=\:\int_{\mathscr{U}}\{\hat{g}_{k}(u)-g_{k0}(u)\}^{2}\mathrm{d}u\\&=\:\int_{\mathscr{U}}\{B^{\mathrm{T}}(u)\hat{\gamma}_{k}-B^{\mathrm{T}}(u)\gamma_{k0}-R_{k}(u)\}^{2}\mathrm{d}u\\&\leqslant2\int_{\mathscr{U}}\{B^{\mathrm{T}}(u)\hat{\gamma}_{k}-B^{\mathrm{T}}(u)\gamma_{k0}\}^{2}\mathrm{d}u+2\int_{\mathscr{U}}R_{k}^{2}(u)\mathrm{d}u\\&=2(\hat{\gamma}_{k}-\gamma_{k0})^{\mathrm{T}}H(\hat{\gamma}_{k}-\gamma_{k0})+2\int_{\mathscr{U}}R_{k}(u)^{2}\mathrm{d}u.\end{aligned}$$
This, together with $\|H\|=O(1)$, can lead to
\begin{equation}
\begin{aligned}
(\hat{\gamma}_{k}-\gamma_{k0})^{\mathrm{T}}H(\hat{\gamma}_{k}-\gamma_{k0})=O_{p}(n^{-2r/(2r+1)}+a_{n}^{2}).\label{tiaojian1}
\end{aligned}
\end{equation}
In addition, it is easy to show that
\begin{equation}
\begin{aligned}
\int_{\mathcal{U}}R_k^2(u)\mathrm{d}u=O_p(n^{-2r/(2r+1)}).\label{tiaojian2}
\end{aligned}
\end{equation}
By using  (\ref{tiaojian1}) and (\ref{tiaojian2}), we complete the proof of (iii).$\hfill\square$\\
\\
$\mathbf{Proof~of~Theorem~2}$\quad (i) By $\lambda _{\max }\to$0, it is easy to show that $a_{n}= 0$ for large $n.$ Then by Theorem 1, it is sufficient to show that, for any $\phi_j$,
$$\|\phi_j-\phi_{j0}\|=O_p(n^{-r/(2r+1)}),\quad j=1,\ldots,s-1,$$
and given small $\epsilon=Cn^{-r/(2r+1)}$, when $n\to\infty$, with probability
approaching to one, it holds that
\begin{equation}
\begin{aligned}
&\frac{\partial Q(\phi,\theta,\gamma)}{\partial\phi_j}>0,\quad 0<\phi_j<\epsilon,\quad j=s,\ldots,p-1,\\
&\frac{\partial Q(\phi,\theta,\gamma)}{\partial\phi_j}<0,\quad -\epsilon<\phi_j<0,\quad j=s,\ldots,p-1.\label{mubiao2}
\end{aligned}
\end{equation}

By a direct calculation, we have
$$\begin{aligned}\frac{\partial Q(\phi,\theta,\gamma)}{\partial\phi_{j}}=~&-\sum_{i=1}^{n}\{Y_{i}-\theta^{\mathrm{T}}U_{i}-W_{i}(\phi)^{\mathrm{T}}\gamma\}\dot{W}_{i}(\phi)^{\mathrm{T}}\gamma e_{\phi_{j}}^{\mathrm{T}}X_{i}+n\dot{p}_{\lambda_{1j}}(|\phi_{j}|)\mathrm{sgn}(\phi_{j})\\=~&-\sum_{i=1}^{n}\{\varepsilon_{i}+R^{\mathrm{T}}(\beta_{0}^{\mathrm{T}}X_{i})Z_{i}+(I_{q}\otimes B(\beta_{0}^{\mathrm{T}}X_{i})\cdot Z_{i})^{\mathrm{T}}(\gamma_{0}-\gamma)\\&+
(I_{q}\otimes[B(\beta_{0}^{\mathrm{T}}X_{i})-B(\beta^{\mathrm{T}}X_{i})]\cdot Z_{i})^{\mathrm{T}}\gamma-(\theta-\theta_0)^{\mathrm{T}}U_{i}\}\{Z_{i}^{\mathrm{T}}\dot{g}(\beta_{0}^{\mathrm{T}}X_{i})\\
&-Z_{i}^{\mathrm{T}}(\dot{g}(\beta_{0}^{\mathrm{T}}X_{i})-I_{q}\otimes\dot{B}^{\mathrm{T}}(\beta_{0}^{\mathrm{T}}X_{i})\gamma_{0})-Z_{i}^{\mathrm{T}}I_{q}\otimes\dot{B}^{\mathrm{T}}(\beta_{0}^{\mathrm{T}}X_{i})(\gamma_{0}-\gamma)\\
&-Z_{i}^{\mathrm{T}}I_{q}\otimes[\dot{B}(\beta_{0}^{\mathrm{T}}X_{i})-\dot{B}(\beta^{\mathrm{T}}X_{i})]^{\mathrm{T}}\gamma\}e_{\phi_{j}}^{\mathrm{T}}X_{i}+n\dot{p}_{\lambda_{1j}}(|\phi_{j}|)\mathrm{sgn}(\phi_{j}),\end{aligned}$$
where $e_{\phi_j}=(-(1-\|\phi\|^2)^{-1/2}\phi_j,0,\ldots,0,1,0,\ldots,0)^\mathrm{T}$
with $(j+1)$th component 1. By conditions (C1), (C2), (\ref{tiaojian}) and Theorem 1, it is easy to show that
$$\frac{\partial Q(\phi,\theta,\gamma)}{\partial\phi_{j}}=n\lambda_{1j}\{\lambda_{1j}^{-1}\dot{p}_{\lambda_{1j}}(|\phi_{j}|)\mathrm{sgn}(\phi_{j})+O_{p}(n^{-r/(2r+1)}\lambda_{1j}^{-1})\}.$$
Since
$$\lim\limits_{n\to\infty}\lim\limits_{\phi_{j}\to0}\inf\lambda_{1j}^{-1}\dot{p}_{\lambda_{1j}}(|\phi_{j}|)>0,\quad\lambda_{1j}n^{r/(2r+1)}\geqslant\lambda_{\min}n^{r/(2r+1)}\to\infty,$$
which means that the sign of the derivative is completely determined by that of $\phi_j.$ Then, we can reach to (\ref{mubiao2}), and the proof of (i) is complete.

 (ii) It can be proved directly by the same arguments in (i).

(iii)  By the similar arguments as in the proof of (i) again,  we have,  with
probability approaching to one, $\hat{\gamma}_k=0,k=d+1,\ldots,q.$ Note that $\hat{g} _k( u) = B^{\mathrm{T} }( u) \hat{\gamma } _k.$  By the fact
\begin{equation}\label{tiaojianB}
\sup_u\|B(u)\|=O(1),
\end{equation}
we can easily  complete the proof.$\hfill\square$\\
\\
$\mathbf{Proof~of~Theorem~3}$\quad By Theorems 1 and 2, we know that, as $n$ $\to\infty$, with probability approaching to one,
$Q(\phi,\theta,\gamma)$ attains the minimal value at $(\hat{\phi}^{*\mathrm{T}},0)^{\mathrm{T}}$, $(\hat{\theta}^{*\mathrm{T}},0)^{\mathrm{T}}$ and $(\hat{\gamma}^{*\mathrm{T}},0)^{\mathrm{T}}$.
Let
$$Q_{1n}(\phi,\theta,\gamma)=\frac{\partial Q(\phi,\theta,\gamma)}{\partial\phi^*}, Q_{2n}(\phi,\theta,\gamma)=\frac{\partial Q(\phi,\theta,\gamma)}{\partial\theta^*}, Q_{3n}(\phi,\theta,\gamma)=\frac{\partial Q(\phi,\theta,\gamma)}{\partial\gamma^*}.$$
Thus, $(\hat{\phi}^{*\mathrm{T}},0)^{\mathrm{T}}$, $(\hat{\theta}^{*\mathrm{T}},0)^{\mathrm{T}}$ and $(\hat{\gamma}^{*\mathrm{T}},0)^{\mathrm{T}}$  satisfy
\begin{align}
&\frac{1}{n}Q_{1n}((\hat{\phi}^{*\mathrm{T}},0)^{\mathrm{T}},(\hat{\theta}^{*\mathrm{T}},0)^{\mathrm{T}},(\hat{\gamma}^{*\mathrm{T}},0)^{\mathrm{T}})\nonumber\\
=~&\frac{-2}{n}\sum_{i=1}^{n}(Y_{i}-\hat{\theta}^{*\mathrm{T}}U_{i}^{*}-W_{i}^{*\mathrm{T}}(\hat{\phi}^{*})\hat{\gamma}^{*})\dot{W}_{i}^{*\mathrm{T}}(\hat{\phi}^{*})\hat{\gamma}^{*}J_{\hat{\phi}^{*}}^{\mathrm{T}}X_{i}^{*}+V_{1}\nonumber\\
=~&0,\lbl{mubiao31}\\
&\frac{1}{n}Q_{2n}((\hat{\phi}^{*\mathrm{T}},0)^{\mathrm{T}},(\hat{\theta}^{*\mathrm{T}},0)^{\mathrm{T}},(\hat{\gamma}^{*\mathrm{T}},0)^{\mathrm{T}})\nonumber\\
=~&\frac{-2}{n}\sum_{i=1}^{n}(Y_{i}-\hat{\theta}^{*\mathrm{T}}U_{i}^{*}-W_{i}^{*\mathrm{T}}(\hat{\phi}^{*})\hat{\gamma}^{*})U_{i}^{*}+V_{2}\nonumber\\
=~&0\lbl{mubiao32}
\end{align}
and
%\end{align}
%\end{equation}
\begin{align}
%\begin{align}
&\frac{1}{n}Q_{3n}((\hat{\phi}^{*\mathrm{T}},0)^{\mathrm{T}},(\hat{\theta}^{*\mathrm{T}},0)^{\mathrm{T}},(\hat{\gamma}^{*\mathrm{T}},0)^{\mathrm{T}})\nonumber\\
=~&\frac{-2}{n}\sum_{i=1}^{n}(Y_{i}-\hat{\theta}^{*\mathrm{T}}U_{i}^{*}-W_{i}^{*\mathrm{T}}(\hat{\phi}^{*})\hat{\gamma}^{*})W_{i}^{*}(\hat{\phi}^{*})+V_{3}\nonumber\\
=~&0,\lbl{mubiao33}
\end{align}
where
\begin{eqnarray*}
\begin{aligned}
&V_1=(\dot{p}_{\lambda_{1,1}}(|\hat{\phi}_1|)\mathrm{sgn}(\hat{\phi}_1),\ldots,\dot{p}_{\lambda_{1,s-1}}(|\hat{\phi}_{s-1}|)\mathrm{sgn}(\hat{\phi}_{s-1}))^\mathrm{T},\\
&V_2=(\dot{p}_{\lambda_{21}}(|\hat{\theta}_1|)\mathrm{sgn}(\hat{\theta}_1),\ldots,\dot{p}_{\lambda_{2w}}(|\hat{\theta}_{w}|)\mathrm{sgn}(\hat{\theta}_{w}))^\mathrm{T},
%{\text{and~}}&V_3=\left(\dot{p}_{\lambda_{31}}(\|\hat{\gamma}_{1}\|_{H})\frac{\hat{\gamma}_{1}^{\mathrm{T}}H}{\|\hat{\gamma}_{1}\|_{H}},\ldots,\dot{p}_{\lambda_{3v}}(\|\hat{\gamma}_{v}\|_{H})\frac{\hat{\gamma}_{v}^{\mathrm{T}}H}{\|\hat{\gamma}_{v}\|_{H}}\right)^{\mathrm{T}}.
\end{aligned}
\end{eqnarray*}
and
\begin{eqnarray*}
\begin{aligned}
V_3=\left(\dot{p}_{\lambda_{31}}(\|\hat{\gamma}_{1}\|_{H})\frac{\hat{\gamma}_{1}^{\mathrm{T}}H}{\|\hat{\gamma}_{1}\|_{H}},\ldots,\dot{p}_{\lambda_{3v}}(\|\hat{\gamma}_{v}\|_{H})\frac{\hat{\gamma}_{v}^{\mathrm{T}}H}{\|\hat{\gamma}_{v}\|_{H}}\right)^{\mathrm{T}}.
\end{aligned}
\end{eqnarray*}

We analyse the above three equations next by the following four steps and try to obtain the relation between the parametric vector $\begin{pmatrix}
\hat{\beta}^{*}-\beta_{0}^{*}\\
\hat{\theta}^{*}-\theta_{0}^{*}
\end{pmatrix}$ and the model error $\varepsilon_i$.

\textbf{Step 1}: We study the components of $V_1$, $V_2$ and $V_3$.

Applying the Taylor expansion to $\dot{p}_{\lambda_{1l}}(|\hat{\phi}_{l}|)$ leads to
$$\dot{p}_{\lambda_{1l}}(|\hat{\phi_l}|)=\dot{p}_{\lambda_{1l}}(|\phi_{l0}|)+\{\ddot{p}_{\lambda_{1l}}(|\phi_{l0}|)+o_{p}(1)\}(\hat{\phi_l}-\phi_{l0}).$$
Furthermore, condition (C5) implies that $\ddot{p}_{\lambda_{1l}}(|\phi_{l0}|)=o_{p}(1)$, and note that $\dot{p}_{\lambda_{1l}}(|\phi_{l0}|)=0$ as $\lambda_{\mathrm{max}}\to0$. Then, it follows by Theorem 1 and Theorem 2 that
$$\dot{p}_{\lambda_{1l}}(|\hat{\phi}_{l}|)\mathrm{sgn}(\hat{\phi}_{l})=o_{p}(\hat{\phi}^{*}-\phi_{0}^{*}).$$
Similarly, we can show
$$\dot{p}_{\lambda_{2h}}(|\hat{\theta}_{h}|)\mathrm{sgn}(\hat{\theta}_{h})=o_{p}(\hat{\theta}^{*}-\theta_{0}^{*}) \text{~~and},$$ $$\dot{p}_{\lambda_{3k}}(\|\hat{\gamma}_{k}\|_{H})\frac{H\hat{\gamma}_{k}}{\|\hat{\gamma}_{k}\|_{H}}=o_{p}(\hat{\gamma}^{*}-\gamma_{0}^{*}).$$
Hence, $V_1=o_{p}(\hat{\phi}^{*}-\phi_{0}^{*})$, $V_2=o_{p}(\hat{\theta}^{*}-\theta_{0}^{*})$ and $V_3=o_{p}(\hat{\gamma}^{*}-\gamma_{0}^{*})$.
~~\\

\textbf{Step 2}:  We obtain the relationship among the parametric vectors
$\hat{\phi}^{*}-\phi_{0}^{*}$, $\hat{\theta}^{*}-\theta_{0}^{*}$ and $\hat{\gamma}^{*}-\gamma_{0}^{*}$  by calculating (\ref{mubiao33}).

A simple calculation yields
$$\begin{aligned}&\frac{1}{n}\sum_{i=1}^{n}(Y_{i}-\hat{\theta}^{*\mathrm{T}}U_{i}^{*}-W_{i}^{*\mathrm{T}}(\hat{\phi}^{*})\hat{\gamma}^{*})W_{i}^{*}(\hat{\phi}^{*})+o_{p}(\hat{\gamma}^{*}-\gamma_{0}^{*})\\
=~&\frac{1}{n}\sum_{i=1}^{n}\left\lbrace{\varepsilon_{i}+R^{\mathrm{T}}(\beta_{i}^{*\mathrm{T}}X_{i}^{*})Z_{i}^{*}-U_{i}^{*\mathrm{T}}(\hat{\theta}^{*}-\theta_0^{*})-W_{i}^{*\mathrm{T}}(\phi_{0}^{*})(\hat{\gamma}^{*}-\gamma_{0}^{*})}\right.\\
&\left.{-[W_{i}^{*}(\hat{\phi}^{*})-W_{i}^{*}(\phi_{0}^{*})]^{\mathrm{T}}\hat{\gamma}^{*}}\right\rbrace\left\{W_{i}^{*}(\phi_{0}^{*})+[W_{i}^{*}(\hat{\phi}^{*})-W_{i}^{*}(\phi_{0}^{*})] \right \}+o_{p}(\hat{\gamma}^{*}-\gamma_{0}^{*}).\\
\end{aligned}$$
Let
$$\begin{aligned}&\Phi_{n}=\frac{1}{n}\sum_{i=1}^{n}W_{i}^{*}(\phi_{0}^{*})W_{i}^{*\mathrm{T}}(\phi_{0}^{*}),\quad\Psi_{n}=\frac{1}{n}\sum_{i=1}^{n}W_{i}^{*}(\phi_{0}^{*})V_{i}^{*\mathrm{T}},\quad G_{n}=\frac{1}{n}\sum_{i=1}^{n}W_{i}^{*}(\phi_{0}^{*})U_{i}^{*\mathrm{T}},\\ &\Lambda_{n}=\frac{1}{n}\sum_{i=1}^{n}W_{i}^{*}(\phi_{0}^{*})(\varepsilon_{i}+R^{\mathrm{T}}(\beta_{0}^{*\mathrm{T}}X_{i}^{*})Z_{i}^{*}).\end{aligned}$$
Then, by condition (C7), (\ref{tiaojianB}) and Theorem 1, we have
\begin{equation}
\begin{aligned}
\hat{\gamma}^*-\gamma_0^*=[\Phi_n+o_p(1)]^{-1}\{\Lambda_n-\Psi_n(\hat{\phi}^*-\phi_0^*)-G_n(\hat{\theta}^*-\theta_0^*)\}.\label{tiaojian3}
\end{aligned}
\end{equation}

\textbf{Step 3}: We obtain the relationship among the parametric vectors
$\hat{\phi}^{*}-\phi_{0}^{*}$, $\hat{\theta}^{*}-\theta_{0}^{*}$ ~and the model error $\varepsilon_i$ by calculating (\ref{mubiao31}).

By substituting (\ref{tiaojian3}) into (\ref{mubiao31}), we can get
\begin{align*}
0=~&\frac{1}{n}\sum_{i=1}^{n}(Y_{i}-\hat{\theta}^{*\mathrm{T}}U_{i}^{*}-W_{i}^{*}(\hat{\phi}^{*})\hat{\gamma}^{*})\dot{W}_{i}^{*}(\hat{\phi}^{*})\hat{\gamma}^{*}J_{\hat{\phi}^{*}}^{\mathrm{T}}X_{i}^{*}+o_{p}(\hat{\phi}^{*}-\phi_{0})\\
=~&\frac{1}{n}\sum_{i=1}^{n}\left\lbrace{\varepsilon_{i}+R^{\mathrm{T}}(\beta_{0}^{\mathrm{T}}X_{i}^{*})Z_{i}^{*}-U_{i}^{*\mathrm{T}}(\hat{\theta}^{*}-\theta_0^{*})-W_{i}^{\mathrm{T}}(\phi_{0}^{*})(\hat{\gamma}^{*}-\gamma_{0}^{*})}\right.\\
&\left.{-[W_{i}^{*}(\hat{\phi}^{*})-W_{i}^{*}(\phi_{0}^{*})]^{\mathrm{T}}\hat{\gamma}^{*}}\right\rbrace\dot{W}_{i}^{\mathrm{T}}(\hat{\phi}^{*})\hat{\gamma}^{*}J_{\hat{\phi}^{*}}^{\mathrm{T}}X_{i}^{*}+o_{p}(\hat{\phi}^{*}-\phi_{0}^{*})\\
%&=\frac{1}{n}\sum_{i=1}^{n}\{\varepsilon_{i}+R^{\mathrm{T}}(\beta_{i}^{\mathrm{T}}X_{i}^{*})Z_{i}^{*}-U_{i}^{*\mathrm{T}}(\hat{\theta}^{*}-\theta_0^{*})-W_{i}^{\mathrm{T}}(\phi_{0}^{*})[\Phi_n+o_p(1)]^{-1}[\Lambda_n-\Psi_n(\hat{\phi}^*-\phi_0^*)\\
%&-G_n(\hat{\theta}^*-\theta_0^*)-[W_{i}^{*}(\hat{\phi}^{*})-W_{i}^{*}(\phi_{0}^{*})]^{\mathrm{T}}\hat{\gamma}^{*}\} W_{i}^{\mathrm{T}}(\hat{\phi}^{*})\hat{\gamma}^{*}J_{\hat{\phi}^{*}}^{\mathrm{T}}X_{i}^{*}+o_{p}(\hat{\phi}^{*}-\phi_{0}^{*})\\
=~&\frac{1}{n}\sum_{i=1}^{n}\left\lbrace{\varepsilon_{i}+R^{\mathrm{T}}(\beta_{0}^{\mathrm{T}}X_{i}^{*})Z_{i}^{*}-U_{i}^{*\mathrm{T}}(\hat{\theta}^{*}-\theta_0^{*})+W_{i}^{*\mathrm{T}}(\phi_{0}^{*})[\Phi_{n}^{-1}+o_{p}(1)]\Psi_{n}(\hat{\phi}^{*}-\phi_{0}^{*})}\right.\\
&+W_{i}^{*\mathrm{T}}(\phi_{0}^{*})[\Phi_{n}^{-1}+o_{p}(1)]G_n(\hat{\theta}^*-\theta_0^*)-W_{i}^{\mathrm{T}}(\phi_{0}^{*})(\Phi_{n}^{-1}+o_{p}(1))\Lambda_{n}\\
&\left.{-[W_{i}^{*}(\hat{\phi}^{*})-W_{i}^{*}(\phi_{0}^{*})]^{\mathrm{T}}\hat{\gamma}^{*}}\right\rbrace\dot{W}_{i}^{\mathrm{T}}(\hat{\phi}^{*})\hat{\gamma}^{*}J_{\hat{\phi}^{*}}^{\mathrm{T}}X_{i}^{*}+o_{p}(\hat{\phi}^{*}-\phi_{0}^{*})\\
=~&\frac{1}{n}\sum_{i=1}^{n}\{\varepsilon_{i}+R^{\mathrm{T}}(\beta_{0}^{\mathrm{T}}X_{i}^{*})Z_{i}^{*}-W_{i}^{*\mathrm{T}}(\phi_{0}^{*})[\Phi_{n}^{-1}+o_{p}(1)]\Lambda_{n}\}\dot{W}_{i}^{\mathrm{T}}(\hat{\phi}^{*})\hat{\gamma}^{*}J_{\hat{\phi}^{*}}^{\mathrm{T}}X_{i}^{*}\\
&+\frac{1}{n}\sum_{i=1}^{n}\{W_{i}^{\mathrm{T}}(\phi_{0}^{*})[\Phi_{n}^{-1}+o_{p}(1)]\Psi_{n}(\hat{\phi}^{*}-\phi_{0}^{*})\}\dot{W}_{i}^{\mathrm{T}}(\hat{\phi}^{*})\hat{\gamma}^{*}J_{\hat{\phi}^{*}}^{\mathrm{T}}X_{i}^{*}\\
&-\frac{1}{n}\sum_{i=1}^{n}[W_{i}^{*}(\hat{\phi}^{*})-W_{i}^{*}(\phi_{0}^{*})]^{\mathrm{T}}\hat{\gamma}^{*}\dot{W}_{i}^{\mathrm{T}}(\hat{\phi}^{*})\hat{\gamma}^{*}J_{\hat{\phi}^{*}}^{\mathrm{T}}X_{i}^{*}\\
&+\frac{1}{n}\sum_{i=1}^{n}\{W_{i}^{\mathrm{T}}(\phi_{0}^{*})[\Phi_{n}^{-1}+o_{p}(1)]G_{n}-U_i^{*{\mathrm{T}}}\}(\hat{\theta}^{*}-\theta_{0}^{*})\dot{W}_{i}^{\mathrm{T}}(\hat{\phi}^{*})\hat{\gamma}^{*}J_{\hat{\phi}^{*}}^{\mathrm{T}}X_{i}^{*}+o_{p}(\hat{\phi}^{*}-\phi_{0}^{*})\\
=:~&J_1+J_2-J_3+J_4+o_p(\hat{\phi}^*-\phi_0^*).\end{align*}

For $J_{1}$, a direct calculation yields
\begin{align*}J_{1}=~&\frac{1}{n}\sum_{i=1}^{n}M_{1i}\boldsymbol{\dot{g}}^{\mathrm{T}}(\beta_{0}^{*\mathrm{T}}X_{i}^{*})Z_{i}^{*}J_{\hat{{\phi}}^{*}}^{\mathrm{T}}X_{i}^{*}\\
&-\frac{1}{n}\sum_{i=1}^{n}M_{1i}[\boldsymbol{\dot{g}}^{\mathrm{T}}(\beta_{0}^{*\mathrm{T}}X_{i}^{*})Z_{i}^{*}-\dot{W}_{i}^{*\mathrm{T}}(\phi_{0}^{*})\gamma_{0}^{*}]J_{\hat{\phi}^{*}}^{\mathrm{T}}X_{i}^{*}\\
&+\frac{1}{n}\sum_{i=1}^{n}M_{1i}\dot{W}_{i}^{*\mathrm{T}}(\phi_{0}^{*})(\hat{\gamma}^{*}-\gamma_{0}^{*})J_{\hat{\phi}^{*}}^{\mathrm{T}}X_{i}^{*}\\
&-\frac{1}{n}\sum_{i=1}^{n}M_{1i}[\dot{W}_{i}^{*}(\phi_{0}^{*})-\dot{W}_{i}^{*}(\hat{\phi}^{*})]^{\mathrm{T}}\hat{\gamma}^{*}J_{\hat{\phi}^{*}}^{\mathrm{T}}X_{i}^{*}\\=:~&J_{11}+J_{12}+J_{13}+J_{14},\end{align*}
where
$M_{1i}=\varepsilon_{i}+R^{\mathrm{T}}(\beta_{0}^{*\mathrm{T}}X_{i}^{*})Z_{i}^{*}-W_{i}^{*\mathrm{T}}(\phi_{0}^{*})[\Phi_{n}^{-1}+o_{p}(1)]\Lambda_{n}.$
Note that
$$\begin{aligned}&\dfrac{1}{n}\sum_{i=1}^{n}\Psi_{n}^{\mathrm{T}}\Phi_{n}^{-1}W_{i}^{*}(\phi_{0}^{*})\{\varepsilon_{i}+R^{\mathrm{T}}(\beta_{0}^{*\mathrm{T}}X_{i}^{*})Z_{i}^{*}-W_{i}^{*\mathrm{T}}(\phi_{0}^{*})\Phi_{n}^{-1}\Lambda_{n}\}=0,\\
&\frac{1}{n}\sum_{i=1}^{n}[V_{i}^{*}-\Psi_{n}^{\mathrm{T}}\Phi_{n}^{-1}W_{i}^{*}(\phi_{0}^{*})]W_{i}^{*\mathrm{T}}(\phi_{0}^{*})=0,\quad J_{\hat{\phi}^{*}}-J_{\phi_{0}^{*}}=O_{p}(\hat{\phi}^{*}-\phi_{0}^{*}).\end{aligned}$$
Then, by condition (C7) and $\|R(u)\|=O(K^{-r})$, we can drive that
$$\begin{aligned}J_{11}
=~&\frac{1}{n}\sum_{i=1}^{n}[V_{i}^{*}-\Psi_{n}^{\mathrm{T}}\Phi_{n}^{-1}W_{i}^{*}(\phi_{0}^{*})]\varepsilon_{i}\\
&+\frac{1}{n}\sum_{i=1}^{n}[V_{i}^{*}-\Psi_{n}^{\mathrm{T}}\Phi_{n}^{-1}W_{i}^{*}(\phi_{0}^{*})]R^{\mathrm{T}}(\beta_{0}^{*\mathrm{T}}X_{i}^{*})Z_{i}^{*}\\
&-\frac{1}{n}\sum_{i=1}^{n}[V_{i}^{*}-\Psi_{n}^{\mathrm{T}}\Phi_{n}^{-1}W_{i}^{*}(\phi_{0}^{*})]W_{i}^{*\mathrm{T}}(\phi_{0}^{*})[\Phi_{n}^{-1}+o_{p}(1)]\Lambda_{n}+o_{p}(\hat{\phi}^{*}-\phi_{0}^{*})\\
=~&\frac{1}{n}\sum_{i=1}^{n}[V_{i}^{*}-\Psi_{n}^{\mathrm{T}}\Phi_{n}^{-1}W_{i}^{*}(\phi_{0}^{*})]\varepsilon_{i}+o_{p}(\hat{\phi}^{*}-\phi_{0}^{*}).\end{aligned}$$
In addition, by (\ref{tiaojian}), it is easy to show that
$$J_{12}=o_p(\hat{\phi}^*-\phi_0^*).$$
Similarly, we can prove that
$$J_{13}=o_p(\hat{\phi}^*-\phi_0^*),\quad J_{14}=o_p(\hat{\phi}^*-\phi_0^*).$$

We now deal with $J_2.$ A simple calculation yields
$$\begin{aligned}
J_{2}=~&\frac{1}{n}\sum_{i=1}^{n}M_{2i}{\boldsymbol{\dot{g}}}^{\mathrm{T}}(\beta_{0}^{*\mathrm{T}}X_{i}^{*})Z_{i}^{*}J_{\hat{\phi}^{*}}^{\mathrm{T}}X_{i}^{*}\\
&-\frac{1}{n}\sum_{i=1}^{n}M_{2i}[\dot{\boldsymbol{g}}^{\mathrm{T}}(\beta_{0}^{*\mathrm{T}}X_{i}^{*})Z_{i}^{*}-W_{i}^{*\mathrm{T}}(\phi_{0}^{*})\gamma_{0}^{*}]J_{\hat{\phi}^{*}}^{\mathrm{T}}X_{i}^{*}\\
&+\frac{1}{n}\sum_{i=1}^{n}M_{2i}\dot{W}_{i}^{*\mathrm{T}}(\phi_{0}^{*})(\dot{\gamma}^{*}-\gamma_{0}^{*})J_{\hat{\phi}^{*}}^{\mathrm{T}}X_{i}^{*}\\
&-\frac{1}{n}\sum_{i=1}^{n}M_{2i}[\dot{W}_{i}^{*}(\phi_{0}^{*})-\dot{W}_{i}^{*}(\hat{\phi}^{*})]^{\mathrm{T}}\hat{\gamma}^{*}J_{\hat{\phi}^{*}}^{\mathrm{T}}X_{i}^{*}\\
=:~&J_{21}+J_{22}+J_{23}+J_{24},\end{aligned}$$
where
$$M_{2i}=W_{i}^{*\mathrm{T}}(\phi_{0}^{*})[\Phi_{n}^{-1}+o_{p}(1)]\Psi_{n}(\hat{\phi}^{*}-\phi_{0}^{*}).$$
By condition (C7), we have
$$J_{21}=\frac{1}{n}\sum_{i=1}^{n}V_{i}^{*}W_{i}^{*\mathrm{T}}(\phi_{0}^{*})\Phi_{n}^{-1}\Psi_{n}(\hat{\phi}^{*}-\phi_{0}^{*})+o_{p}(\hat{\phi}^{*}-\phi_{0}^{*}).$$
Similar arguments  of $J_{12}$  can lead to
$$J_{22}=o_p(\hat{\phi}^*-\phi_0^*),\quad J_{23}=o_p(\hat{\phi}^*-\phi_0^*),\quad J_{24}=o_p(\hat{\phi}^*-\phi_0^*).$$

We now consider $J_3.$ By (\ref{tiaojian}), Theorem 1, and the Taylor expansion
technique, we have
$$\begin{aligned}J_{3}
=~&\frac{1}{n}\sum_{i=1}^{n}[\dot{W}_{i}^{*\mathrm{T}}(\phi_{0}^{*})\hat{\gamma}^{*}X_{i}^{*\mathrm{T}}J_{\phi_{0}^*}(\hat{\phi}^{*}-\hat{\phi}_{0}^{*})+o_{p}(\hat{\phi}^{*}-\phi_{0}^{*})]\dot{W}_{i}^{*\mathrm{T}}(\hat{\phi}^{*})\hat{\gamma}^{*}J_{\hat{\phi}^{*}}^{\mathrm{T}}X_{i}^{*}\\
=~&\frac{1}{n}\sum_{i=1}^{n}[\boldsymbol{\dot{g}}^{\mathrm{T}}(\beta_{0}^{\mathrm{T}}X_{i}^{*})Z_{i}^{*}X_{i}^{\mathrm{T}}J_{\phi_{0}^*}(\hat{\phi}^{*}-\phi_{0}^{*})+o_{p}(\hat{\phi}^{*}-\phi_{0}^{*})]\dot{W}_{i}^{*\mathrm{T}}(\hat{\phi}^{*})\hat{\gamma}^{*}J_{\hat{\phi}^{*}}^{\mathrm{T}}X_{i}^{*}\\
=~&\frac{1}{n}\sum_{i=1}^{n}V_{i}^{*}V_{i}^{*}(\hat{\phi}^{*}-\phi_{0}^{*})-\frac{1}{n}\sum_{i=1}^{n}V_{i}^{*}(\hat{\phi}^{*}-\phi_{0}^{*})[\boldsymbol{\dot{g}}^{\mathrm{T}}(\beta_{0}^{*\mathrm{T}}X_{i}^{*})Z_{i}^{*}\\
&-\dot{W}_{i}^{*\mathrm{T}}(\phi_{0}^{*})\gamma_{0}^{*}
]J_{\hat{\phi}^{*}}^{\mathrm{T}}X_{i}^{*}+\frac{1}{n}\sum_{i=1}^{n}V_{i}^{*}(\hat{\phi}^{*}-\phi_{0}^{*}))\dot{W}_{i}^{*\mathrm{T}}(\phi_{0}^{*})(\hat{\gamma}^{*}-\gamma_{0}^{*})J_{\hat{\phi}^{*}}^{\mathrm{T}}X_{i}^{*}\\
&-\frac{1}{n}\sum_{i=1}^{n}V_{i}^{*\mathrm{T}}(\hat{\phi}^{*}-\phi_{0}^{*})[\dot{W}_{i}^{*}(\phi_{0}^{*})-\dot{W}_{i}^{*}(\hat{\phi}^{*})]^{\mathrm{T}}\hat{\gamma}^{*}J_{\hat{\phi}^{*}}^{\mathrm{T}}X_{i}^{*}+o_{p}(\hat{\phi}^{*}-\phi_{0}^{*})\\
=~&\frac{1}{n}\sum_{i=1}^{n}V_{i}V_{i}^{*\mathrm{T}}(\hat{\phi}^{*}-\phi_{0}^{*})+o_{p}(\hat{\phi}^{*}-\phi_{0}^{*}).
\end{aligned}$$

%$$\begin{aligned}&\sqrt{n}\left(\hat{\phi}^{*}-\phi_{0}^{*}\right)\\=~&\left\{\frac{1}{n}\sum_{i=1}^{n}[V_{i}^{*}V_{i}^{*\mathrm{T}}-\Psi_{n}^{\mathrm{T}}\Phi_{n}^{-1}\Psi_{n}]+o_{p}(1)\right\}^{-1}\frac{1}{\sqrt{n}}\sum_{i=1}^{n}[V_{i}^{*}-\Psi_{n}^{\mathrm{T}}\Phi_{n}^{-1}W_{i}^{*}(\phi_{0}^{*})]\varepsilon_{i}.\end{aligned}$$
For the term $J_{4}$,
\begin{equation*}
\begin{aligned}J_{4}
=~&\frac{1}{n}\sum_{i=1}^{n}\{W_{i}^{\mathrm{T}}(\phi_{0}^{*})[\Phi_{n}^{-1}+o_{p}(1)]G_{n}-U_i^{*{\mathrm{T}}}](\hat{\theta}^{*}-\theta_{0}^{*})\}W_{i}^{\mathrm{T}}(\hat{\phi}^{*})\hat{\gamma}^{*}J_{\hat{\phi}^{*}}^{\mathrm{T}}X_{i}^{*}
\end{aligned}
\end{equation*}
\begin{equation*}
\begin{aligned}
J_{4}=~&\frac{1}{n}\sum_{i=1}^{n}M_{4i}\dot{\boldsymbol{g}}^{\mathrm{T}}(\beta_{0}^{*\mathrm{T}}X_{i}^{*})Z_{i}^{*}J_{\hat{\phi}^{*}}^{\mathrm{T}}X_{i}^{*}\\
&-\frac{1}{n}\sum_{i=1}^{n}M_{4i}[\dot{\boldsymbol{g}}^{\mathrm{T}}(\beta_{0}^{*\mathrm{T}}X_{i}^{*})Z_{i}^{*}-W_{i}^{*\mathrm{T}}(\phi_{0}^{*})\gamma_{0}^{*}]J_{\hat{\phi}^{*}}^{\mathrm{T}}X_{i}^{*}\\
&+\frac{1}{n}\sum_{i=1}^{n}M_{4i}\dot{W}_{i}^{*\mathrm{T}}(\phi_{0}^{*})(\hat{\gamma}^{*}-\gamma_{0}^{*})J_{\hat{\phi}^{*}}^{\mathrm{T}}X_{i}^{*}\\
&-\frac{1}{n}\sum_{i=1}^{n}M_{4i}[\dot{W}_{i}^{*}(\phi_{0}^{*})-\dot{W}_{i}^{*}(\hat{\phi}^{*})]^{\mathrm{T}}\hat{\gamma}^{*}J_{\hat{\phi}^{*}}^{\mathrm{T}}X_{i}^{*}\\
=:~&J_{41}+J_{42}+J_{43}+J_{44},
\end{aligned}
\end{equation*}
where
$M_{4i}=\{W_{i}^{\mathrm{T}}(\phi_{0}^{*})[\Phi_{n}^{-1}+o_{p}(1)]G_{n}-U_i^{*{\mathrm{T}}}\}(\hat{\theta}^{*}-\theta_{0}^{*}).$
$$\begin{aligned}
J_{41}=~&\frac{1}{n}\sum_{i=1}^{n}V_{i}^{*}W_{i}^{*\mathrm{T}}(\phi_{0}^{*})\Phi_{n}^{-1}G_{n}(\hat{\theta}^{*}-\theta_{0}^{*})-\frac{1}{n}\sum_{i=1}^{n}V_{i}^{*}U_{i}^{*\mathrm{T}}(\hat{\theta}^{*}-\theta_{0}^{*})+o_{p}(\hat{\theta}^{*}-\theta_{0}^{*})\\
%=~&(\Psi_{n}^{\mathrm{T}}\Phi_{n}^{-1}G_{n}-\frac{1}{n}\sum_{i=1}^{n}V_{i}^{*}U_{i}^{*\mathrm{T}})(\hat{\theta}^{*}-\theta_{0}^{*})+o_p(\hat{\theta}^{*}-\theta_{0}^{*}).
\end{aligned}$$
By (\ref{tiaojian}) and (C7), it is easy to show that
$$J_{42}=o_p(\hat{\theta}^*-\theta_0^*).$$
Similarly, we can prove that
$$J_{43}=o_p(\hat{\theta}^*-\theta_0^*),\quad J_{44}=o_p(\hat{\theta}^*-\theta_0^*).$$

From the above arguments, we get
\begin{align}
&\frac{1}{n}\sum_{i=1}^{n}V_{i}^{*}[V_{i}^{*\mathrm{T}}-W_{i}^{*\mathrm{T}}(\phi_{0}^{*})\Phi_{n}^{-1}\Psi_{n}](\hat{\phi}^{*}-\phi_{0}^{*})-\frac{1}{n}\sum_{i=1}^{n}V_{i}^{*}W_{i}^{*\mathrm{T}}(\phi_{0}^{*})\Phi_{n}^{-1}G_{n}(\hat{\theta}^{*}-\theta_{0}^{*})\nonumber\\
&+\frac{1}{n}\sum_{i=1}^{n}V_{i}^{*}U_{i}^{*\mathrm{T}}(\hat{\theta}^{*}-\theta_{0}^{*})
+o_{p}(\hat{\theta}^{*}-\theta_{0}^{*})+o_{p}(\hat{\phi}^{*}-\phi_{0}^{*})\nonumber\\
=~&\frac{1}{n}\sum_{i=1}^{n}\left\{V_{i}^{*}V_{i}^{*\mathrm{T}}-\Psi_{n}^{\mathrm{T}}\Phi_{n}^{-1}\Psi_{n}\right\}(\hat{\phi}^{*}-\phi_{0}^{*})+\frac{1}{n}\sum_{i=1}^{n}\left\{V_{i}^{*}U_{i}^{*\mathrm{T}}-\Psi_{n}^{\mathrm{T}}\Phi_{n}^{-1}G_{n}\right\}(\hat{\theta}^{*}-\theta_{0}^{*})\nonumber\\
&+o_p(\hat{\theta}^{*}-\theta_{0}^{*})+o_{p}(\hat{\phi}^{*}-\phi_{0}^{*})\nonumber\\
=~&\frac{1}{n}\sum_{i=1}^{n}[V_{i}^{*}-\Psi_{n}^{\mathrm{T}}\Phi_{n}^{-1}W_{i}^{*}(\phi_{0}^{*})]\varepsilon_{i}.\label{jielun1}
\end{align}
%ÁíÆðÒ»¶Î

\textbf{Step 4}: We  obtain the other relation among the parametric vectors
$\hat{\phi}^{*}-\phi_{0}^{*}$, $\hat{\theta}^{*}-\theta_{0}^{*}$ ~and the model error $\varepsilon_i$ by calculating (\ref{mubiao32}).

Substituting (\ref{tiaojian3}) into (\ref{mubiao32}), we can get
\begin{align}
0=~&\frac{1}{n}\sum_{i=1}^{n}(Y_{i}-\hat{\theta}^{*\mathrm{T}}U_{i}^{*}-W_{i}^{*}(\hat{\phi}^{*}))U_{i}^*+o_{p}(\hat{\theta}^{*}-\theta_{0})\nonumber\\
=~&\frac{1}{n}\sum_{i=1}^{n}\{\varepsilon_{i}+R^{\mathrm{T}}(\beta_{0}^{\mathrm{T}}X_{i}^{*})Z_{i}^{*}-U_{i}^{*\mathrm{T}}(\hat{\theta}^{*}-\theta_0^{*})-W_{i}^{\mathrm{T}}(\phi_{0}^{*})(\hat{\gamma}^{*}-\gamma_{0}^{*})\nonumber\\&-[W_{i}^{*}(\hat{\phi}^{*})-W_{i}^{*}(\phi_{0}^{*})]^{\mathrm{T}}\hat{\gamma}^{*}\}U_{i}^*+o_{p}(\hat{\theta}^{*}-\theta_{0}^{*})\nonumber\\
%&=\frac{1}{n}\sum_{i=1}^{n}\{\varepsilon_{i}+R^{\mathrm{T}}(\beta_{0}^{\mathrm{T}}X_{i}^{*})Z_{i}^{*}-U_{i}^{*\mathrm{T}}(\hat{\theta}^{*}-\theta_0^{*})-W_{i}^{\mathrm{T}}(\phi_{0}^{*})[\Phi_n+o_p(1)]^{-1}[\Lambda_n-\Psi_n(\hat{\phi}^*-\phi_0^*)\\
%&-G_n(\hat{\theta}^*-\theta_0^*)-[W_{i}^{*}(\hat{\phi}^{*})-W_{i}^{*}(\phi_{0}^{*})]^{\mathrm{T}}\hat{\gamma}^{*}\} W_{i}^{\mathrm{T}}(\hat{\phi}^{*})\hat{\gamma}^{*}J_{\hat{\phi}^{*}}^{\mathrm{T}}X_{i}^{*}+o_{p}(\hat{\phi}^{*}-\phi_{0}^{*})\\
=~&\frac{1}{n}\sum_{i=1}^{n}\{\varepsilon_{i}+R^{\mathrm{T}}(\beta_{0}^{\mathrm{T}}X_{i}^{*})Z_{i}^{*}-U_{i}^{*\mathrm{T}}(\hat{\theta}^{*}-\theta_0^{*})+W_{i}^{*\mathrm{T}}(\phi_{0}^{*})[\Phi_{n}^{-1}+o_{p}(1)]\Psi_{n}(\hat{\phi}^{*}-\phi_{0}^{*})\nonumber\\
&+W_{i}^{*\mathrm{T}}(\phi_{0}^{*})[\Phi_{n}^{-1}+o_{p}(1)]G_n(\hat{\theta}^*-\theta_0^*)-W_{i}^{\mathrm{T}}(\phi_{0}^{*})(\Phi_{n}^{-1}+o_{p}(1))\Lambda_{n}\nonumber\\&-[W_{i}^{*}(\hat{\phi}^{*})-W_{i}^{*}(\phi_{0}^{*})]^{\mathrm{T}}\hat{\gamma}^{*}\}U^{*}+o_{p}(\hat{\theta}^{*}-\theta_{0}^{*})\nonumber\\
=~&\frac{1}{n}\sum_{i=1}^{n}\{\varepsilon_{i}+R^{\mathrm{T}}(\beta_{0}^{\mathrm{T}}X_{i}^{*})Z_{i}^{*}-W_{i}^{*\mathrm{T}}(\phi_{0}^{*})[\Phi_{n}^{-1}+o_{p}(1)]\Lambda_{n}]U_i^{*}\nonumber\\
&+\frac{1}{n}\sum_{i=1}^{n}\{W_{i}^{\mathrm{T}}(\phi_{0}^{*})[\Phi_{n}^{-1}+o_{p}(1)]\Psi_{n}(\hat{\phi}^{*}-\phi_{0}^{*})\}U_{i}^{*}
-\frac{1}{n}\sum_{i=1}^{n}[W_{i}^{*}(\hat{\phi}^{*})-W_{i}^{*}(\phi_{0}^{*})]^{\mathrm{T}}U_{i}^{*}\nonumber\\
&+\frac{1}{n}\sum_{i=1}^{n}\{W_{i}^{\mathrm{T}}(\phi_{0}^{*})[\Phi_{n}^{-1}+o_{p}(1)]G_{n}-U_i^{*{\mathrm{T}}}\}(\hat{\theta}^{*}-\theta_{0}^{*})U_{i}^{*}+o_{p}(\hat{\theta}^{*}-\theta_{0}^{*})\nonumber\\
=:~&J_1^{'}+J_2^{'}-J_3^{'}+J_4^{'}+o_p(\hat{\theta}^{*}-\theta_{0}^{*}).\nonumber
\end{align}

For $J_{1}^{'}$, a direct calculation yields
\begin{align*}J_{1}^{'}=~&\frac{1}{n}\sum_{i=1}^{n}M_{1i}U_i^*\end{align*}
where
$M_{1i}=\varepsilon_{i}+R^{\mathrm{T}}(\beta_{0}^{*\mathrm{T}}X_{i}^{*})Z_{i}^{*}-W_{i}^{*\mathrm{T}}(\phi_{0}^{*})[\Phi_{n}^{-1}+o_{p}(1)]\Lambda_{n}.$
Note that
$$\begin{aligned}&\dfrac{1}{n}\sum_{i=1}^{n}G_{n}^{\mathrm{T}}\Phi_{n}^{-1}W_{i}^{*}(\phi_{0}^{*})\{\varepsilon_{i}+R^{\mathrm{T}}(\beta_{0}^{*\mathrm{T}}X_{i}^{*})Z_{i}^{*}-W_{i}^{*\mathrm{T}}(\phi_{0}^{*})\Phi_{n}^{-1}\Lambda_{n}\}=0,\\
&\frac{1}{n}\sum_{i=1}^{n}[U_{i}^{*}-G_{n}^{\mathrm{T}}\Phi_{n}^{-1}W_{i}^{*}(\phi_{0}^{*})]W_{i}^{*\mathrm{T}}(\phi_{0}^{*})=0.\end{aligned}$$
By condition (C7) and $\|R(u)\|=O(K^{-r})$, we can drive that
$$\begin{aligned}J_{1}^{'}=~&\frac{1}{n}\sum_{i=1}^{n}[U_{i}^{*}-G_{n}^{\mathrm{T}}\Phi_{n}^{-1}W_{i}^{*}(\phi_{0}^{*})]\varepsilon_{i}\\&+\frac{1}{n}\sum_{i=1}^{n}[U_{i}^{*}-G_{n}^{\mathrm{T}}\Phi_{n}^{-1}W_{i}^{*}(\phi_{0}^{*})]R^{\mathrm{T}}(\beta_{0}^{*\mathrm{T}}X_{i}^{*})Z_{i}^{*}\\&-\frac{1}{n}\sum_{i=1}^{n}[U_{i}^{*}-G_{n}^{\mathrm{T}}\Phi_{n}^{-1}W_{i}^{*}(\phi_{0}^{*})]W_{i}^{*\mathrm{T}}(\phi_{0}^{*})[\Phi_{n}^{-1}+o_{p}(1)]\Lambda_{n}+o_{p}(\hat{\phi}^{*}-\phi_{0}^{*})\\=~&\frac{1}{n}\sum_{i=1}^{n}[U_{i}^{*}-G_{n}^{\mathrm{T}}\Phi_{n}^{-1}W_{i}^{*}(\phi_{0}^{*})]\varepsilon_{i}.\end{aligned}$$

For $J_{2}^{'}$, a simple calculation yields
%$$\begin{aligned}J_{2}=\dfrac{1}{n}\sum_{i=1}^{n}G_{n}^{\mathrm{T}}\Phi_{n}^{-1}\Psi_n(\hat{\phi}^{*}-\phi_{0}^{*})+o_{p}(\hat{\phi}^{*}-\phi_{0}^{*}),\end{aligned}$$
$$\begin{aligned}J_{2}^{'}=\dfrac{1}{n}\sum_{i=1}^{n}G_{n}^{\mathrm{T}}\Phi_{n}^{-1}\Psi_n(\hat{\phi}^{*}-\phi_{0}^{*})+o_{p}(\hat{\phi}^{*}-\phi_{0}^{*}).\end{aligned}$$

Similarly, we also can get
$$\begin{aligned}J_{3}^{'}
=~&\frac{1}{n}\sum_{i=1}^{n}[\dot{W}_{i}^{*\mathrm{T}}(\phi_{0}^{*})\hat{\gamma}^{*}X_{i}^{*\mathrm{T}}J_{\phi_{0}^*}(\hat{\phi}^{*}-\hat{\phi}_{0}^{*})+o_{p}(\hat{\phi}^{*}-\phi_{0}^{*})]U_i^*\\
=~&\frac{1}{n}\sum_{i=1}^{n}[\boldsymbol{\dot{g}}^{\mathrm{T}}(\beta_{0}^{\mathrm{T}}X_{i}^{*})Z_{i}^{*}X_{i}^{\mathrm{T}}J_{\phi_{0}^*}(\hat{\phi}^{*}-\phi_{0}^{*})+o_{p}(\hat{\phi}^{*}-\phi_{0}^{*})]U_i^*\\
=~&\frac{1}{n}\sum_{i=1}^{n}U_{i}^{*}V_{i}^{*\mathrm{T}}(\hat{\phi}^{*}-\phi_{0}^{*})+o_{p}(\hat{\phi}^{*}-\phi_{0}^{*})
\end{aligned}$$
$$\begin{aligned}J_{4}^{'}=(G_{n}^{\mathrm{T}}\Phi_{n}^{-1}G_{n}-\frac{1}{n}\sum_{i=1}^{n}U_{i}^{*}U_{i}^{*\mathrm{T}})(\hat{\theta}^{*}-\theta_{0}^{*})+o_{p}(\hat{\theta}^{*}-\theta_{0}^{*}),\end{aligned}$$

From the above arguments, we can obtain
\begin{align}
&\frac{1}{n}\sum_{i=1}^{n}\left\{U_{i}^{*}V_{i}^{*\mathrm{T}}-G_{n}^{\mathrm{T}}\Phi_{n}^{-1}\Psi_{n}\right\}(\hat{\phi}^{*}-\phi_{0}^{*})+\frac{1}{n}\sum_{i=1}^{n}\left\{U_{i}^{*}U_{i}^{*\mathrm{T}}-G_{n}^{\mathrm{T}}\Phi_{n}^{-1}G_{n}\right\}(\hat{\theta}^{*}-\theta_{0}^{*})\nonumber\\
&+o_p(\hat{\theta}^{*}-\theta_{0}^{*})+o_{p}(\hat{\phi}^{*}-\phi_{0}^{*})\nonumber\\
=~&\frac{1}{n}\sum_{i=1}^{n}[U_{i}^{*}-G_{n}^{\mathrm{T}}\Phi_{n}^{-1}W_{i}^{*}(\phi_{0}^{*})]\varepsilon_{i}.\label{jielun2}
\end{align}

Combining (\ref{jielun1}) with (\ref{jielun2}), we can obtain
%\begin{aligned}
%&\begin{pmatrix}  \frac{1}{n}\sum_{i=1}^{n}[V_{i}^{*}V_{i}^{*\mathrm{T}}-\Psi_{n}^{\mathrm{T}}\Phi_{n}^{-1}\Psi_{n}]+o_p(1)
%&\frac{1}{n}\sum_{i=1}^{n}[V_{i}^{*}U_{i}^{*\mathrm{T}}-\Psi_{n}^{\mathrm{T}}\Phi_{n}^{-1}G_{n}]+o_p(1) \\  \frac{1}{n}\sum_{i=1}^{n}[U_{i}^{*}V_{i}^{*\mathrm{T}}-G_{n}^{\mathrm{T}}\Phi_{n}^{-1}\Psi_{n}]+o_p(1)
%&\frac{1}{n}\sum_{i=1}^{n}[U_{i}^{*}U_{i}^{*\mathrm{T}}-G_{n}^{\mathrm{T}}\Phi_{n}^{-1}G_{n}]+o_p(1)
%\end{pmatrix}
%\begin{pmatrix}
%\hat{\phi}^{*}-\phi_{0}^{*}\\
%\hat{\theta}^{*}-\theta_{0}^{*}
%\end{pmatrix}\\
%=~&\begin{pmatrix}\frac{1}{n}\sum_{i=1}^{n}[V_{i}^{*}-\Psi_{n}^{\mathrm{T}}\Phi_{n}^{-1}W_{i}^{*}(\phi_{0}^{*})]\\
%\frac{1}{n}\sum_{i=1}^{n}[U_{i}^{*}-G_{n}^{\mathrm{T}}\Phi_{n}^{-1}W_{i}^{*}(\phi_{0}^{*})]
%\end{pmatrix}\varepsilon_{i}
%\end{aligned}$$
$$
\begin{aligned}
&\sqrt{n}
\begin{pmatrix}
\hat{\phi}^{*}-\phi_{0}^{*}\\
\hat{\theta}^{*}-\theta_{0}^{*}
\end{pmatrix}\\
=~&{\begin{pmatrix}
 \frac{1}{n}\sum_{i=1}^{n}[V_{i}^{*}V_{i}^{*\mathrm{T}}-\Psi_{n}^{\mathrm{T}}\Phi_{n}^{-1}\Psi_{n}]+o_p(1)
&\frac{1}{n}\sum_{i=1}^{n}[V_{i}^{*}U_{i}^{*\mathrm{T}}-\Psi_{n}^{\mathrm{T}}\Phi_{n}^{-1}G_{n}]+o_p(1) \\  \frac{1}{n}\sum_{i=1}^{n}[U_{i}^{*}V_{i}^{*\mathrm{T}}-G_{n}^{\mathrm{T}}\Phi_{n}^{-1}\Psi_{n}]+o_p(1)
&\frac{1}{n}\sum_{i=1}^{n}[U_{i}^{*}U_{i}^{*\mathrm{T}}-G_{n}^{\mathrm{T}}\Phi_{n}^{-1}G_{n}]+o_p(1)
\end{pmatrix}}^{-1}\\
&\cdot\frac{1}{\sqrt{n}}\begin{pmatrix}\sum_{i=1}^{n}[V_{i}^{*}-\Psi_{n}^{\mathrm{T}}\Phi_{n}^{-1}W_{i}^{*}(\phi_{0}^{*})]\\
\sum_{i=1}^{n}[U_{i}^{*}-G_{n}^{\mathrm{T}}\Phi_{n}^{-1}W_{i}^{*}(\phi_{0}^{*})]
\end{pmatrix}\varepsilon_{i}.
\end{aligned}$$

It now follows from (\ref{2.3}) that
$$\hat{\beta}^{*}-\beta_{0}^{*}=J_{\phi_{0}^{*}}(\hat{\phi}^{*}-\phi_{0}^{*})+O_{p}(n^{-1}).$$
Thus, we have
$$
\begin{aligned}
&\sqrt{n}
\begin{pmatrix}
\hat{\beta}^{*}-\beta_{0}^{*}\\
\hat{\theta}^{*}-\theta_{0}^{*}
\end{pmatrix}\\
=~&\begin{pmatrix}
  J_{\phi_{0}^{*}}& 0\\
 0 &I_w
\end{pmatrix}{\begin{pmatrix}
 \frac{1}{n}\sum_{i=1}^{n}V_{i}^{*}V_{i}^{*\mathrm{T}}-\Psi_{n}^{\mathrm{T}}\Phi_{n}^{-1}\Psi_{n}
&\frac{1}{n}\sum_{i=1}^{n}V_{i}^{*}U_{i}^{*\mathrm{T}}-\Psi_{n}^{\mathrm{T}}\Phi_{n}^{-1}G_{n} \\  \frac{1}{n}\sum_{i=1}^{n}U_{i}^{*}V_{i}^{*\mathrm{T}}-G_{n}^{\mathrm{T}}\Phi_{n}^{-1}\Psi_{n}
&\frac{1}{n}\sum_{i=1}^{n}U_{i}^{*}U_{i}^{*\mathrm{T}}-G_{n}^{\mathrm{T}}\Phi_{n}^{-1}G_{n}
\end{pmatrix}}^{-1}\\
&\cdot\frac{1}{\sqrt{n}}\begin{pmatrix}\sum_{i=1}^{n}[V_{i}^{*}-\Psi_{n}^{\mathrm{T}}\Phi_{n}^{-1}W_{i}^{*}(\phi_{0}^{*})]\\
\sum_{i=1}^{n}[U_{i}^{*}-G_{n}^{\mathrm{T}}\Phi_{n}^{-1}W_{i}^{*}(\phi_{0}^{*})]
\end{pmatrix}\varepsilon_{i}+o_p(1).
\end{aligned}$$

Combining the central limit theorem with Slutsky's
theorem, we can easily obtain by Lemma \ref{yl2} that
%combining with Lemma 2, by the central limit theorem and Slutsky's theorem, we obtain
$$\begin{aligned}
&\sqrt{n}
\begin{pmatrix}
\hat{\beta}^{*}-\beta_{0}^{*}\\
\hat{\theta}^{*}-\theta_{0}^{*}
\end{pmatrix}\stackrel{\mathcal{L}}{\longrightarrow}N(0,\sigma^{2}\widetilde{J}_{\phi_{0}^{*}}\Sigma^{-1}\widetilde{J}_{\phi_{0}^{*}}^{\mathrm{T}}),\end{aligned}$$
where $\widetilde{J}_{\phi_{0}^{*}}=\begin{pmatrix}
  J_{\phi_{0}^{*}}& 0\\
 0 &I_w
\end{pmatrix}.$
This completes the proof of Theorem 3.$\hfill\square$

%\section*{Acknowledgements}
%The authors would like to thank two anonymous referees for many valuable comments and suggestions that have led to improvements in the paper.

\end{document}